\documentclass[times]{cnmauth}
\usepackage{moreverb}

\usepackage[colorlinks,bookmarksopen,bookmarksnumbered,citecolor=red,urlcolor=red]{hyperref}
\usepackage{subfigure}

\newcommand\BibTeX{{\rmfamily B\kern-.05em \textsc{i\kern-.025em b}\kern-.08em
T\kern-.1667em\lower.7ex\hbox{E}\kern-.125emX}}

\def\volumeyear{}

\newcommand\BigR{\mathbb{R}}

\newcommand\BigE{\mathbb{E}}

\newcommand\Acal{\mathcal{A}}

\newcommand\Dcal{\mathcal{D}}

\newcommand\Pcal{\mathcal{P}}

\newcommand\Ucal{\mathcal{U}}
\newcommand\Mcal{\mathcal{M}}

\newcommand\bu{\boldsymbol{u}}

\newcommand\bxi{\boldsymbol{\xi}}

\newcommand\bU{\boldsymbol{U}}
\newcommand\bH{\boldsymbol{H}}

\newcommand\KL{Karhunen-Lo\`eve }
\newcommand\cov{\textrm{cov}}

\begin{document}

\runningheads{A. Brault}{UQ of pulse wave propagation in a vascular network}

\title{Uncertainty quantification of inflow boundary condition and proximal arterial stiffness coupled effect on\\ pulse wave propagation in a vascular network}

\author{A. Brault $^\dagger$, L. Dumas$^\ddagger$ and D. Lucor\corrauth}

\address{
$^\dagger$
Institut de Math\'ematiques de Toulouse, Universit\'e Paul Sabatier, 118, Route de Narbonne, 31062 Toulouse cedex 4, France
\\
$^\ddagger$
Laboratoire de Math\'ematiques de Versailles, UVSQ, CNRS, Universit\'e Paris-Saclay, 78035 Versailles, France
\\
$^{*}$
LIMSI, CNRS, Universit\'e Paris-Saclay, Campus universitaire b\^at 508, Rue John von Neumann, F-91405 Orsay cedex, France\\
}

\corraddr{LIMSI, CNRS, Universit\'e Paris-Saclay, Campus universitaire b\^at 508, Rue John von Neumann, F-91405 Orsay cedex, France,  \texttt{didier.lucor@limsi.fr}}

\begin{abstract}
This work aims at quantifying the effect of inherent uncertainties from cardiac output on the sensitivity of a human compliant arterial network response based on stochastic simulations of a reduced-order pulse wave propagation model.
A simple pulsatile output form is utilized to reproduce the most relevant cardiac features with a minimum number of  parameters associated with left ventricle dynamics. Another source of critical uncertainty is the spatial heterogeneity of the aortic compliance which plays a key role in the propagation and damping of pulse waves generated at each cardiac cycle. A continuous representation of the aortic stiffness in the form of a generic random field of prescribed spatial correlation is then considered. 
Resorting to a stochastic sparse pseudospectral method, we investigate the spatial sensitivity of the pulse pressure and waves reflection magnitude with respect to the different model uncertainties.
Results indicate that uncertainties related to the shape and magnitude of the prescribed inlet flow in the proximal aorta can lead to potent variation of both the mean value and standard deviation of blood flow velocity and pressure dynamics due to the interaction of different wave propagation and reflection features. These results have potential physiological and pathological implications. They will provide some guidance in clinical data acquisition and future coupling of arterial pulse wave propagation reduced-order model with more complex beating heart models.
\end{abstract}


\maketitle

%

\section{Introduction}
Human cardiovascular system and {its function and structure are} primarily affected by aging, eating habits, life style and other  risk factors, inducing increased stiffness which is associated to {cardiovascular} and cerebral morbi-mortality. 
Arterial blood pressure and flow waves are generated by the heart contraction and
its interaction with a continuous network of variable distensible vessels (the arteries, microcirculation and veins). The arteries, that carry blood from the heart, form a predominantly tree-like structure, called the arterial tree. The (large) arteries distend to accommodate to the blood volume change and generate {\em pulse} waves that induce continuous variations of blood pressure and velocity from central to distal locations of the arterial network. The mechanism is very complex as pulse waveforms and their reflections depend for instance on geometry, elastic properties, local pulse wave speeds, velocity profiles and boundary conditions.

Several recent publications \cite{Sankaran2011,Chen2013,Eck_2015} have demonstrated the interest of uncertainty quantification (UQ) of the circulatory system since many aleatory and epistemic uncertainties emerge due to its complexity, diversity and multiscale nature. Unfortunately, these new techniques are nowadays still out of reach for full 3D {\em patient-specific} computational fluid dynamics (CFD) simulations involving detailed time-dependent hemodynamics. For now, it must rely on reduced-order models (ROM) instead, owing to the computational burden related to the curse of dimensionality induced by the parametric lack of knowledge and uncertainty. Reduced-order pulse wave propagation models,  
based on the integrated forms of the Navier-Stokes equations of  continuity and momentum and a description of the arterial wall compliance, are good compromise between accuracy and computational cost, and have been validated against {\em in vivo} data, e.g. \cite{Reymond_2011,Bollache}. These distributed models are computationally efficient to simulate wave propagation phenomena. Because of their flexibility in predictions of pressure and flow in all  main arteries, they have the potential to assist practicians in the diagnosis and treatment of cardiovascular pathologies, in particular if they are properly calibrated with subject-specific data. There have been two different classes of approach for the UQ of such systems. In the first one, the ROM is treated as a generic model that would represent the ``average" response of a group of individuals \cite{Chen2013,Eck_2015}. In this case, the effect of uncertainty is accounted for global sensitivity analysis of several input parameters. This allows the identification of the model parameters that influence the output the most. It provides better general understanding and links those trends to general potential physiological and pathological implications. Nevertheless, the system under consideration is sometimes too generic and the uncertainty modeling too simple and idealized, e.g. disregarding the correlations and physical constraints among parameters.\\
On the other hand, some studies are concerned with a more specific model for a given patient for which one holds clinical data/measurements \cite{Leguy2011,Auricchio2015}. In this case, the parameters are typically inferred from experiments that provide only indirect observations: the outputs of these experiments are influenced by the parameters which can then be estimated by solving an {\em inverse problem}. When measurement error is included in the observations and uncertainty considered in the \emph{forward model} describing the effect of the sought-after parameters onto the outputs, it results in a stochastic inverse problem whose solution is a probabilistic description of the parameters, in contrast with a point-wise estimate as in inverse problems \cite{Auricchio2015,OlmBayes2016}.\\
An intermediate tackle of the problem is to consider a given (subject or patient) model for which we assume we know (most of) the parameters but where the uncertainty resides in the operating conditions (e.g. boundary conditions) that affects the model in response to a range of different external sollicitations (e.g. physical activities).

In this work, we represent the human cardiovascular system as a simple ROM in the form of a distributed one-dimensional arterial network wave propagation fluid-structure interaction model. The model encompasses the 55 main arteries of the human body. This type of model is very sensitive to the variability of the mechanical geometric and elastic properties of the arterial walls, as well as the outflow boundary conditions (BC) that are imposed to the truncated and simplified computational network \cite{FormaggiaBCs2002,VignonClementel20063776,Grinberg_ABE2008}. While recent works raise the question of the relevance of imposed inflow boundary conditions -- i.e. idealized vs. personalized velocity profiles -- in computational hemodynamics \cite{Marzo2011,Morbiducci_2013,Mynard2013870,XiangBCs2014}, very little has been validated for ROM of the arterial tree \cite{Franke_PhD2003}.\\
In terms of uncertainty quantification, some work have been carried out to study the effect of the terminal outflow boundary conditions \cite{Chen2013} but nothing exists, to our knowledge, that reports the effect of inlet boundary conditions waveforms. This is the first contribution that we propose in this paper.
In particular the inflow BC (at the proximal aorta level) that mimics the cardiac output is subject to uncertainties arising from the shape and magnitude of the prescribed flow rate. This inflow is in fact function of some physiological parameters associated with left ventricle dynamics including heart rate ($hr$), or equivalently  cardiac cycle period $T=\frac{1}{hr}$, mean cardiac output $\overline{Q}$, and an estimation of the peak-to-mean flow ratio. These parameters are patient-specific and time-dependent and are therefore difficult to calibrate. The parameter $hr$ at rest is for instance a significant risk factor for cardiovascular and non-cardiovascular death \cite{Albaladejo2001}.\\
We propose to investigate the impact of the uncertainty inherent to the cardiac output by considering $T$, $\overline{Q}$ and peak flow rate $Q_{\textrm{max}}$ as independent uniformly distributed random variables (RV). The quantities of interest (QoI) that we monitor include pulse pressure, arterial distensibility and the pulse reflection magnitude, which are known to be determined in part by the ejection pattern of the left ventricle.\\
It is known that the compliance of the aorta, that connects the left ventricle of the heart to the systemic circulation, also plays a key role in the propagation and damping of the pressure waves generated at each cardiac cycle. The exact mechanism by which arterial stiffening, as a result of ageing or arterial disease, leads or contributes to pulse wave velocity increase, primary hypertension \cite{Laurent01052001} and left ventricular hypertrophy \cite{ORourke_2007} remains a controversial topic. Moreover, stiffness properties of the aortic walls are not well known and subject to large (spatial) variabilities among individuals. Indeed, the elastic properties of conduit arteries vary along the arterial tree; with more elastic proximal arteries and stiffer distal arteries. It is therefore important to account for some uncertainty in the aortic stiffness of our model while preserving a physiological sound heterogeneity in the streamwise spatial distribution of the arterial wall properties. 

Our second contribution in this paper, is
to generalize the representation of the aortic stiffness in the form of generic random {\em fields} with imposed correlation lengths. The main improvement with this approach is that each realization of the random field is continuous in space. In particular, there is no random and unphysical jump across the interface between two contiguous arterial segments in the considered domain. Moreover, it is easy to explore different scenarios depending on the level of spatial correlation considered to the aortic stiffness along its span.\\
In this work, we make use of  {\em non-intrusive} uncertainty propagation techniques and the system stochastic responses are approximated thanks to orthogonal pseudospectral polynomial surrogates of the random input parameters. Those continuous surrogates are constructed from the problem solution sampling that we wish as sparse and efficient as possible. To this end, we rely on the use of sparse grids \cite{SG_book} together with a careful  adjustment of the span of our approximation functional space.\\
In the following, we first introduce the deterministic reduced-order model of the arterial circulation, together with its initial and boundary conditions, that is later carefully validated in the results section. Parametric stochastic modeling and subsequent stochastic response approximation and {\em global} sensitivity analysis are presented in the next section. Finally, the numerical method and posterior analysis are first applied  to the case of the uncertainty quantification of stochastic inflow boundary condition on pulse waves response in a human arterial network. This setup is then revisited with the incorporation of a stochastic modeling of the aortic stiffness in order to study the outcome of the coupling of these two sources of uncertainties.

\section{Methodology}
In the following, we introduce the  simplified blood flow model of a human arterial network and how we couple it with a stochastic approximation of the time-dependent inflow BC that is used to account for the pulsatile cardiac output uncertainty delivered to the system. Later on, this stochastic modeling is enriched by considering that the aortic portion of the network naturally exhibits biological spatial variation in arterial stiffness. In particular, we monitor the effect of the coupled interaction between these two sources of uncertainties -- i.e. upstream cardiac output BC and -- proximal aortic stiffness characterized as a random field, on the pressure waves.
\begin{figure}
\centering
\includegraphics[width=14.0cm]{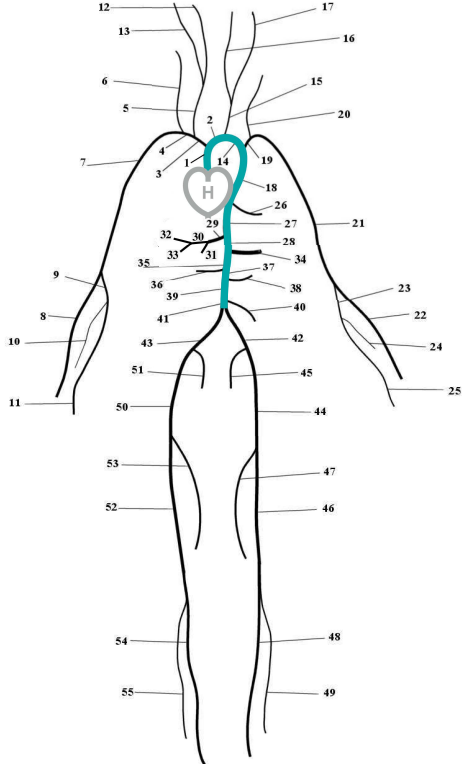}
\caption{Diagrammatic representation of the connectivity and numbering of the 55-largest arteries human network, described in details in ~\cite{Sherwin_2003a}. Uncertainty modeling is introduced for the cardiac output (H) first and is then combined with the aorta stiffness material properties (i.e. for arteries $\mathrm{A}_{\{1,2,14,18,27,28,35,37,39,41\}}$ represented in cyan color).}
 \label{fig:network}
\end{figure}

\subsection{Reduced-order model of the arterial circulation}
\subsubsection{Governing equations}
We consider  the main systemic arterial tree based on $n=55$ arteries  described in the work of \cite{Sherwin_2003a}: we have a network of thin, deformable, {and axisymmetric} arterial segments filled with blood, taken as an incompressible Newtonian fluid ~\cite{Sherwin_2003a}. The formulation, for each arterial segment, based on the conservation of mass and momentum laws and on Young-Laplace equation, is:
\begin{eqnarray}
  \frac{\partial A}{\partial t} + \frac{\partial Au}{\partial x} & = & 0 \nonumber \\
  \frac{\partial u}{\partial t} + u \frac{\partial u}{\partial x} & = & - \frac{1}{\rho} \frac{\partial p}{\partial x} + \frac{f}{\rho \,A} \nonumber \\
  p & = & p_{0} + \beta \big ( \sqrt{A}-\sqrt{A_0} \big ),
\label{eq:govEq}
\end{eqnarray}
where $t$ denotes time, $x \in \Dcal$ is the axial coordinate along the arterial centerline, $A(x,t)$ is the circular cross-sectional area of the lumen, $u(x,t)$ and $p(x,t)$ are average velocity and internal pressure, respectively, $\rho$ is blood density. The term $f$ is the friction force per unit length and is related to the velocity profile through $\displaystyle f = -2\mu \pi \frac{\alpha}{\alpha -1}u$, where $\mu$ is blood dynamic viscosity and $\alpha\in [0,1[$ is a correction factor accounting for the nonlinear integration of radial velocities in each cross-section ~\cite{Smith}. The underscript $_0$ denotes quantities at rest: {$p_{0}$} is the external pressure outside the arterial segment, $A_0$ is the luminal cross-section at pressure $p_{0}$; $A_0$ and $h_0$ are taken {\em constant} along each segment. The {$\beta$} parameter is a measure of the arterial wall stiffness related to its mechanical behaviour: 
\begin{equation}
\beta(x)=\sqrt{\pi} \,h_0 E(x)/{(1-\nu^2)A_0},
\end{equation}
where $h_0$ is the reference arterial wall thickness, $E$ is the Young's modulus that may be variable along the domain and $\nu=1/2$ is Poisson's ratio. In our study, {$\beta$} has therefore a spatial distribution that does not depend on time and may be subject to substantial biological variability among individuals. It is also related to the local distensibility $\displaystyle D=\frac{1}{A} \frac{dA}{dp}$ or pulse wave velocity $c$ through the following relation: 
\begin{eqnarray}
\beta(x) = \frac{2}{D\sqrt{A}} = 2\rho \, \frac{c^2(x,t)}{\sqrt{A(x,t)}}= 2\rho \, \frac{c_0^2(x)}{\sqrt{A_0}},
\label{eq:beta_pwv_dist}
\end{eqnarray}
which  implicitly requires that: $ \displaystyle c_0^2(x) \equiv \frac{\sqrt{\pi} \,h_0 E(x)}{{2\rho(1-\nu^2)}\sqrt{A_0}}$.

Note that $c$ and $c_0$ increase with increasing elastic modulus and wall thickness, and decreasing luminal area.\\
Moreover, we may express: $\displaystyle \frac{\partial p}{\partial x} = \frac{\partial p}{\partial A}\frac{\partial A}{\partial x} + \frac{\partial p}{\partial \beta}\frac{d\beta}{dx}= \frac{\partial p}{\partial A}\frac{\partial A}{\partial x} + \frac{\partial p}{\partial c_0}\frac{dc_0}{dx}$ and write the system of equations in the following non-conservative form:
\begin{equation}
\label{eq:syshyp1}
\frac{\partial \bU}{\partial t}+\bH(\bU) \, \frac{\partial \bU}{\partial x}=
\left[
\begin{array} {c}
A\\
u
\end{array}
\right]_t 
+
\left[
\begin{array} {cc}
u & A\\
c^2/A & u
\end{array}
\right]
\left[
\begin{array} {c}
A\\
u
\end{array}
\right]_x 
=
{\boldsymbol{S} },
\end{equation}
where:

\begin{equation}
{\boldsymbol{S}=\left[
\begin{array} {c}
0\\
\frac{f}{\rho A} - 4 \left (\sqrt{\frac{A}{A_0}} -1 \right ) c_0  \frac{dc_0}{dx} \end{array}
\right].}
\end{equation}

\subsubsection{Deterministic boundary conditions}
\label{subsubsec:outlet}
Perfusion of the peripheral micro-circulation must be accounted in the form of simplified time-dependent lumped models. In the following, 
RCR-type zero-dimensional (0D) outlet boundary conditions are used for the network studied in this paper \cite{Alastruey_2008}. The outlet lumped parameter model accounts for effects of distal arteries located beyond the {simulated} segment. We assume that this missing part is modeled as a tube of length $l$. In this case there exists a friction term governed by Poiseuille's flow, which is a better approximation in small arteries than the flatter velocity profile used in large arteries \cite{Olufsen_2000}. Starting from the linearized version of system (\ref{eq:govEq}) \cite{Womersley_1957} around its reference state (written in $(Q,p)$ variables), integrating along $l$ and {using the assumption} that blood inertia forces are now negligible, we end up with the Windkessel model \cite{Westerhof_2009}:
\begin{equation}
\frac{dp_{in}}{dt} + \frac{p_{in}}{RC} = \frac{Q_{in}}{C} +\frac{p_{v}}{RC},
\label{eq:eq0D}
\end{equation}
where
$\displaystyle R = 8\pi \mu l / A_0^2$ and $\displaystyle C = A_0l / \rho c_0^2$
are the {total} arterial resistance and wall compliance, respectively, $Q_{in}$ is provided as an output from the arterial network reduced-order model, 
$p_v$ is the outflow venous pressure and $p_{in}$ is an
approximation of the mean arterial pressure.\\
\begin{figure}[htbp]
\begin{center}
\includegraphics[width=0.35\textwidth]{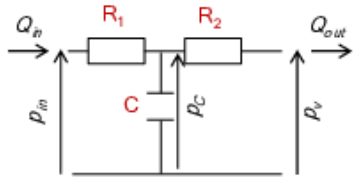}
\caption{\label{fig:RCR} {Electric circuit analogy with lumped 0D $RCR$ Windkessel model.}}
\end{center}
\end{figure}
The three-element RCR version of this approximation, cf. Figure~\ref{fig:RCR}, call{ed here} ${R_1CR_2}$ model, is composed of a resistance alone ${R_{1}}$ coupled with a ${R_{2}C}$ model.
Resistance $\displaystyle {R_{1}}=\rho\, c_0 / A_0$ is added to eliminate {non-physiological} reflected pulse wave oscillations, 
while the ${R_{2}C}$ model is governed by:
\begin{equation}
\frac{dp_{C}}{dt} = \frac{1}{C}\left (Q_{in} +\frac{p_{v}-p_C}{R_2} \right ),
\label{eq:p_C_1}
\end{equation}
{where $p_C$ is pressure accross $C$.}
The ${R_1CR_2}$ model is coupled to our simulated 1D arterial segment at the outlet through {the upwinded flow rate}:
\begin{equation}
Q^u = A^uu^u = \frac{p(A^u) - p_C}{R_{1}},
\label{eq:Auuu_coupl01D}
\end{equation}
which, using invariance of characteristic variables across the outlet{,} leads to the following nonlinear equation:
\begin{eqnarray}
R_{1}\left (u_l + 4\sqrt{\frac{\beta}{2\rho}}A_l^{1/4} \right )A^u - 4R_{1}\sqrt{\frac{\beta}{2\rho}}(A^u)^{5/4} \nonumber \\
- \beta\, \big ((A^u)^{1/2} - A_0^{1/2}\big ) -p_{0} +p_v = 0.
\end{eqnarray}
This equation is solved using a Newton-Raphson method with $A^u=A_l$ as initial guess providing $u^u$ from Eq. (\ref{eq:Auuu_coupl01D}), which serves imposing fluid velocity BC $u_r$ (assuming $A_l=A_r$) such that $u_r = 2u^u - u_l$. 
Finally, $p_C$ is determined at every time step $n$ from Eq. (\ref{eq:p_C_1}). Structural uncertainty modeling will be considered later on, for the aortic portion of the network only. This part of the system not being connected to outflow BCs, will preclude the RCR models from turning to stochastic entities.\\
At the arteries bifurcations, conservation of mass and continuity of the total pressure between the outlet of the parent vessel and the inlet of the daughter vessels are used to close the nonlinear system of equations and avoid dissipation of energy in the arterial branching \cite{Formaggia2003}.\\

\subsubsection{Numerical solver}
The method of characteristics \cite{Alastruey_2012} (not presented here) is a useful analysis prior to the formulation of the numerical scheme.
In this work, we rely on a discontinuous Galerkin (DG) method with a spectral~$/hp$ spatial discretization ~\cite{Karniadakis_2005}, {which} is a very efficient scheme for high-order discretization of convection-dominated flows. The numerical scheme propagates waves at different frequencies without suffering from excessive dispersion and diffusion errors. Detailed mathematical formulation and analysis about the hyperbolic system and the numerical scheme may be found for instance in \cite{Alastruey_2012,Bollache}.

We use an explicit second-order Adams-Bashforth scheme for time integration and inner products are evaluated by numerical integration. Our choice of time step is dictated by the CFL condition and the number of Gauss-Legendre points follows the $3/2$ over-integration rule ~\cite{Karniadakis_2005}. In the case of polynomial approximation within each cell of order $p=3$, pulse pressure (PP) convergence is reached for an average mesh resolution of 13 cells per unit meter. 
~\\
This model requires an accurate calculation of the upwinded fluxes across inter-elemental boundaries. We {tested} several formulations to solve the Riemann problem arising at each element interfaces, including exact and approximate resolution{s} ~\cite{Toro_2009}. {Since} no significant differences were {found regarding} results accuracy under physiological conditions, the Riemann solver of Roe, {which is} numerically more efficient, {was used}.\\
{T}ime-dependent boundary conditions (BC) at the inlet of the modelled segment {to model flow rate from the heart (cf. Section \ref{subsubsec:intlet})} and outlet (to model the missing distal arterial domain, cf. Section \ref{subsubsec:outlet}) {were specified}. 

For both cases, we face a Riemann problem at the interface between the beginning or end of the arterial segment and a virtual location outside the segment. At time $t$, an interface separates two {initial} states for area and velocity {$(A_l,u_l)$} on the left and {$(A_r,u_r)$} on the right, resulting at $t+{\Delta}t$ in two {\it upwinded} states. While supposing rarefaction waves, characteristic variables are constant between the initial and upwinded states.
Based on these assumptions, the setup is quite straightforward referring to \cite{Alastruey_2006} for numerical implementation details.\\
\subsubsection{Wave separation calculation}
\label{subsubsec:wave_separation}
A wave separation technique \cite{Parker_1990} has been deployed for a better interpretation of the main pressure and flow waves. This technique, that accounts for the nonlinear effects due to the convective acceleration and the relation between pressure and arterial change in cross-section, allows the separation between the forward $(\cdot)_f$ and backward $(\cdot)_b$ wave contributions. In the following, we will be interested by an estimation of the reflexion magnitude: i.e. the $p_b/p_f$ ratio at different network locations.
The ratio estimation requires proper integration constants to use in the integration of the difference data to obtain the measured data. In our case, the integration constant for both forward and backward components is defined as half of the mean value of the measured signal over one cardiac cycle.

\subsubsection{Initial condition}
\begin{figure}
\centering
\includegraphics[width=9cm]{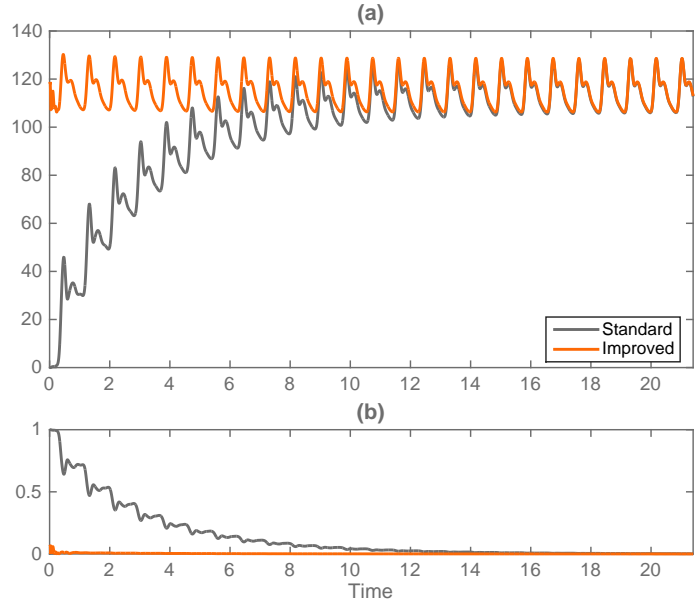}
\caption{Example of temporal buildup of (a): simulated pressure waveforms in the thoracic aorta for different types of cross-sectional area initialization. Standard initialization corresponds to $A(x,t=0)=A_0$, and improved version to Eq.~(\ref{eq:A_init}) and (b): corresponding relative errors with respect to the asymptotic periodic stable solution.}
 \label{fig:initial_condition}
\end{figure}

Time recording of pressure, velocity and flow rate signals show that several heartbeat cycles (e.g. 10 to 20) are in general necessary for the simulation to reach a stable periodic state. In particular, pulse waveforms become periodic when all the blood volume flow accommodated by the network (equipped with its RCR BCs) during systole is released during diastole \cite{Alastruey_2006}. This comes at a considerable cost when one performs a statistical analysis of the system, because full simulations are carried out numerous times for different input parameters. Shortening each simulation will result in substantial computational savings.\\
We propose to initialize the network with arterial cross-sectional areas values that are more representative of the asymptotically periodic pulse pressure. We derive our approximation from the 0D Windkessel model version of the whole systemic arterial tree that we integrate over one cardiac period. In this case, we express it from the pressure law as:
\begin{equation}
A(x,t=0)=\left ( \frac{\overline{p_w}-p_0}{\beta(x)}+ \sqrt{A_0} \right )^2,
\label{eq:A_init}
\end{equation}
where the integrated Windkessel pressure is defined as $\overline{p_w}=R_T \times \overline{Q}+p_v$, with $\overline{Q}\equiv \overline{Q_{out}} = \overline{Q_{\mathrm{inlet}}}$ and $R_T$ is the total resistance of the system\footnote{The total resistance $R_T$ must include resistance of each artery $\mathrm{A}_{\{i\}}$ of the discretized network due to friction: $R_{f_i}=2(\zeta+2)\pi \mu \,l_{\mathrm{A}_{\{i\}}}/A_{0_i}^2$, $\zeta=(2-\alpha)/(\alpha-1)$ ($\zeta=9$ in our case) as well as the equivalent resistance $R_{RCR_i}$ of each RCR outflow BC. The total resistance is found by reducing the different series and parallel combinations step-by-step to end up with a single equivalent resistance for the network.}. In our case, this trick provides good results and the simulations reach a periodic state much faster. Only about 3  cardiac cycles are necessary in order to do so for all arteries in the network, cf. Figure~\ref{fig:initial_condition} for an example.

\subsection{Stochastic modeling} 
\label{subsec:stoch_mod}
\subsubsection{Inflow boundary conditions}
\label{subsubsec:intlet}

A periodic functional form for the inflow BC into the ascending aorta is used to represent the pulsatile cardiac output based on heart rate and stroke volume and capture most relevant features observed in measured physical data, cf. work of Stevens {\em et al.} \cite{Stevens_MB_2003}. It is completely specified by a minimal number of three standard input parameters associated with left ventricle dynamics, namely cardiac cycle period $T$, mean cardiac output $\overline{Q}$ (sometimes named CO in the literature), and maximal flow rate  $Q_{\textrm{max}}$, where $Q_{\textrm{max}}\simeq 6\, \overline Q$ for healthy resting adults. Heart rate (affected by age, fitness levels, hormones and autonomic innervation) and stroke volume (affected by heart size, gender, contractility, duration of contraction, preload and afterload) are the primary determinants of cardiac output. In fact, mean arterial pressure (MAP) should be directly proportional to $\overline{Q}$. The chosen model function describes the time-dependent arterial flow rate (in ml/s) at the inlet:

\begin{eqnarray}
Q_{\textrm{inlet}} (t,\omega,n,\phi)=\Acal\times \overline{Q}\sin^n(\omega t) \cos(\omega t - \phi),
\label{eq:Qinlet}
\end{eqnarray}
where $\omega=\frac{\pi}{T}$, $\phi \in ]0,\pi/2]$ is  the phase angle  and the term $\Acal$ comes from a normalizing factor: $\Acal=\sqrt{\pi}\,\Gamma(\frac{3+n}{2})/ \big(\Gamma(\frac{1+n}{2})\sin(\phi)\big)$.

A typical value of $n=13$ is selected for this study. It represents the right systolic span (about 1/3 of the cardiac cycle) over the cardiac period range. Denoting $Q_{\textrm{max}}$ the maximum flow rate value (at $t=t^*$), the phase angle $\phi$
 is then a solution of the implicit relation:
 \begin{equation}
Q_{\textrm{inlet}} (t^*,\omega,n,\phi )=Q_{\textrm{max}}.
\label{eq:Qinlet2}
\end{equation}
 
 In practice, once a triplet $(T,\overline{Q},Q_{\textrm{max}})$ is chosen, the previous equations are numerically solved for $\phi$ and the expression for time-dependent $Q_{\textrm{inlet}}$ is then readily available, thanks to Eq. (\ref{eq:Qinlet}).\\
 
\begin{figure}
\centering
\includegraphics[width=9cm]{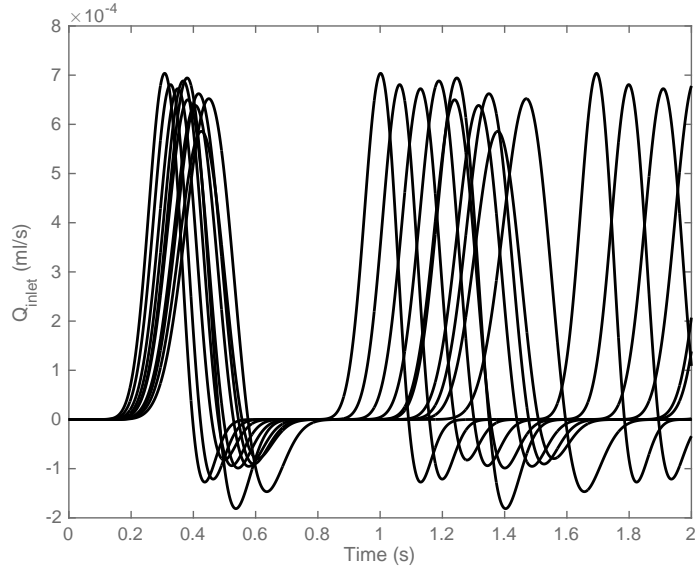}
\caption{Exemple of several random realizations of aortic ascending inflow rates $Q_{\textrm{inlet}}$ versus time for various values of $(T,\overline{Q},Q_{\textrm{max}})$.}
 \label{fig:Q_profiles}
\end{figure}

 Due to our lack of knowledge and thorough statistical analysis, the physical quantities of this triplet are considered in the following as independent, uniformly distributed $\sim \Ucal_{[\mu\pm \sigma]}$ random variables (RV)  within reasonable physiological data ranges, cf. low-order statistics: mean and standard deviation (std) in Table (\ref{tab:stats}). They are defined via zero-mean, unit variance uniform random variables $\bxi$. Despite moderate uncertainty ranges, the imposed flow rate exhibits significant variability in time, cf. Figure~\ref{fig:Q_profiles}. For fixed heart beat and fixed mean cardiac output $\overline{Q}$ , an increase in $Q_{\textrm{max}}$ is adjusted by a corresponding increase in the reverse flow magnitude; while for a fixed  $Q_{\textrm{max}}$, an increase in mean cardiac output leads to a decrease in the reverse flow magnitude.

Equations (\ref{eq:Qinlet}-\ref{eq:Qinlet2}) seem to suggest that $Q_{\textrm{max}}$ depend on the choice of $T$ and $\bar{Q}$.
However, we have checked, from a numerical point of vue, that the three parameters could be chosen independently as long as their values remain within reasonable physiological conditions.

\subsubsection{Aortic random stiffness field}
\label{subsubsec:aort_rand}
The exact mechanism by which arterial stiffening, as a result of ageing or arterial disease, leads to increase in PP and systolic hypertension remains a controversial topic. It is known that the compliance of the proximal aorta plays a key role in the propagation and damping of the pressure waves generated at each cardiac cycle. Indeed, about half of the total systemic arterial compliance is located upstream of the proximal thoracic aorta.
Nevertheless, the stiffness properties of the aortic walls are not well known and subject to large (spatial) variabilities among individuals. Recent studies have considered the effect of local proximal aortic stiffening by increasing the characteristic impedance by opposition to a global stiffening of the entire arterial tree. This was empirically and uniformly done by scaling down the distensibility of the upper aorta. This is particularly delicate because the elastic properties of conduit arteries vary along the arterial tree; with more elastic proximal arteries and stiffer distal arteries. This heterogeneity in the aortic stiffness has important physiological and pathophysiological consequences.\\
Some UQ studies exist that have confirmed that the inherent aortic stiffness uncertainty (either due to biological variability or age-related) had a strong influence on the timing and amplitude of the forward and backward propagating waves, e.g. \cite{Chen2013,Eck_2015}. 
For instance, Leguy and coworkers have shown in a global sensitivity analysis of a wave propagation model for arm arteries that the Young's modulus was the most influential parameter in terms of standard cardiovascular outputs such as pulse pressure, brachial to radial pulse pressure amplification index and some waves transit times across the arterial tree \cite{Leguy2011}.
In these studies, the aortic stiffening variability is often treated in the form of an additive or multiplicative independent RV assigned to the stiffness of each or a group of arteries.  
Apart from the consequent added computational cost due to the handling of numerous RV, the main limitation of these approaches is that it neglects the spatial stiffness correlation among neighboring segments of the aorta. Indeed, when material properties are piecewise randomly sampled in the distributed domain, it is likely that the arterial stiffness will exhibit unphysical jumps across bifurcations that will affect the ``stiffness gradient", known to be one of the leading causes of wave reflections. Moreover, stiffer distal (relative to proximal) locations in the uncertain aorta must be physiologically supported, at least in a {\em statistical} sense (i.e. mean profile should grow and arguably std should decay along the aorta). This distribution is not always guarantee based on the realizations of a set of independent RV.\\
In this work, we propose to generalize the representation of the aortic stiffness in the form of generic random fields (or stochastic processes). Thanks to some regularity assumptions, it is possible to discretize this random field with a finite {\em minimum} number of random variables (that we call again $\bxi$), depending on the choice of {\em correlation} kernel and length for the process.\\
The \KL (KL) representation \cite{Karhunen1946,Loeve77} of a stochastic process with finite variance is a type of model reduction similar to an infinite Fourier-type series based on a spectral decomposition of its correlation kernel in an orthogonal coordinate system. It is closely related to similar decomposition techniques such as Proper Orthogonal Decomposition (POD), Principal Component Analysis (PCA) or Singular Value Decomposition (SVD). 
Its $n-$truncated version is optimal in the sense that there is no other {\em linear}
representation with $n$ terms that has a smaller mean square quadratic error. Here we consider that the stiffness of the aorta may be represented as a stationary random field with a spatial correlation defined by a Gaussian kernel. We introduce the covariance model through the pulse wave velocity at rest, $c_0$:
\begin{eqnarray}
\cov_{c_0}: (x_{1},x_{2}) \in \Dcal \times \Dcal \longmapsto \cov_{c_0}(x_{1},x_{2}) & = &
\BigE[({c_0}(x_{1},\omega) - \mu_{c_0}(x_{1}) )({c_0}(x_{2},\omega) - \mu_{c_0}(x_{2}) )] \nonumber \\
& = & \sigma_{c_0}(x) \, e^{\frac{-(x_2-x_1)^2}{2\,C_{l}^2}}.
\end{eqnarray}
\noindent Because $\cov_{c_0}$ is continuous and bounded on $\Dcal \times \Dcal$, the covariance operator is {\em real, symmetric} and {\em positive-definite} and has a
countable sequence of eigenpairs $(\lambda_i,r_i) \in \BigR^+ \times L^2(\Dcal)$, satisfying an eigenvalue problem described by a Fredholm equation of the second-kind:
\begin{equation}
\label{eq:Fredholm}
\int_{\mathcal{D}} \cov_R({x}_{1}, {x}_{2})
r_i({x}_{2}) d{x}_{2} = \lambda_i
r_i({x}_{1}) \; \; \textrm{with} \; \;
\int_{\mathcal{D}} r_i({x}) r_j({x}) d{x}=\delta_{ij},
\end{equation}
and the eigenfunctions form a complete Hilbertian orthogonal basis of $L^2(\Dcal)$.

\begin{figure}
\centering
\subfigure[\KL eigen-spectrum representation: eigenvalues $\lambda_n$ (a), and first eigenfunctions $r_n(x)$ (b).]{\includegraphics[width=7.0cm]{./brault_figure5a.eps}}\hspace{0.1in}
\subfigure[Aortic stiffness random field: realizations along aorta span (cyan color); mean profile (bold solid line) \& std envelop (dashed lines); discrete values from \cite{Stergiopulos_1992} (crosses); $C_l$ is a third of the aorta span $l_{\mathrm{aorta}}$.]{\includegraphics[width=7.0cm]{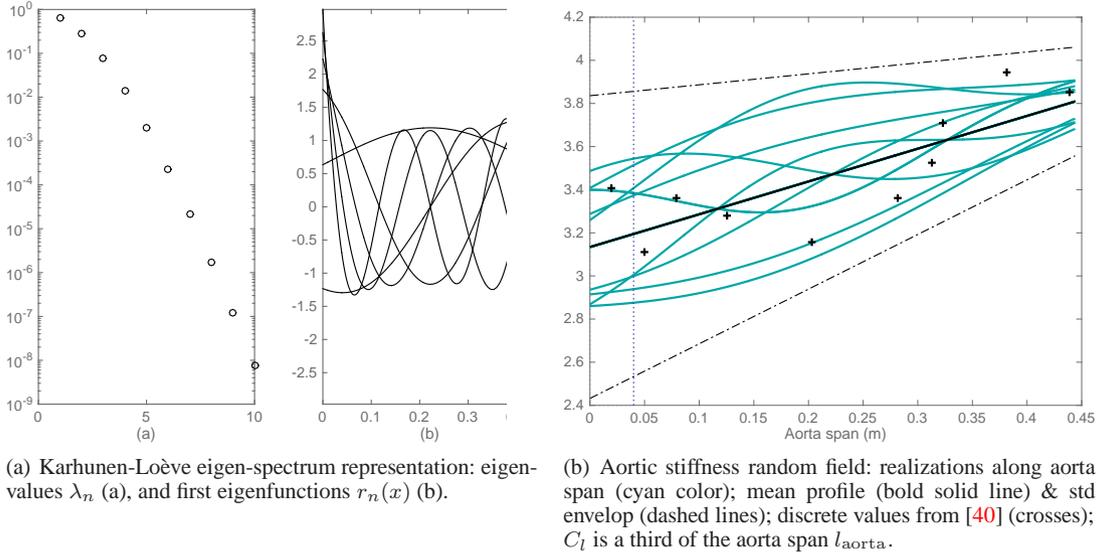}}
\caption{Example of several random realizations of the aorta stiffness random field: here, depicted as pulse wave velocities at rest $c_0$ (b) obtained from a \KL modal representation (a).}
\label{fig:KL_decomposition}
\end{figure}

\noindent The truncated \KL representation of $c_0(x, \bxi)$ is:
\begin{equation}
\label{eq:KL_expansion} c_0(x, \bxi) \approx c_{0_0}(x) + \sum_{i=1}^{n} c_{0_i}(x,\xi_i) = \mu_{c_0}(x) +
\sigma_{c_0}(x) \sum_{i=1}^{n} \sqrt{\lambda_i}\, r_i(x)\, \xi_i,
\end{equation}
where $\bxi$ are mutually {\em uncorrelated} RV with zero mean and unit variance and the $\lambda_i$ are real, ordered, positive eigenvalues. The truncation error decreases monotonically with the number of terms in the expansion. The convergence is inversely proportional to the {\em
correlation length} and depends on the regularity of the covariance kernel. For our case, we have considered that the random array $\bxi$ has components that are uniformly distributed and that a typical correlation length of the random process would have to be quite large in order to mimic the stiffness distribution provided by the literature (cf. cross markers in  Figure~\ref{fig:KL_decomposition}-(c))
To this end, we have chosen a correlation length $C_l \equiv l_{\mathrm{aorta}}/3$, which is well represented with the first $n=3$ modes, cf. Figure~\ref{fig:KL_decomposition}-(a) where one may see eigenvalues and eigenfunctions obtained for this choice of correlation kernel.\\
The main improvement relative to the literature is that the random variability is continuous in space. In particular, there is no random and unphysical jump across the interface between two contiguous arterial segments in the considered domain. Figure~\ref{fig:KL_decomposition}-(b) shows some random realizations of the arterial stiffness along the aortic span, together with its low-order statistics. 

\subsection{Stochastic response approximation}
Now that we have defined, thanks to a set of random variables $\bxi$, the potential uncertain input parameters of our wave propagation reduced-order model, we wish to quantify their effects on some system responses of interest, that we will generically name $\bu$. The pulse wave propagation model $\Mcal$ nonlinearly maps these stochastic variables to the response $\bu = \Mcal({\bxi})$.
In this work, we make use of  {\em non-intrusive} uncertainty propagation techniques: i.e. the sampling of model solutions is done using our solver for the deterministic problem as a black box. We then build a surrogate model that approximates the exact response as accurately as possible based on the smallest number of \emph{observations} or \emph{samples}. Once the samples are acquired, the construction and interrogation of the surrogate model itself is computationally efficient so that the predictive capability of the metamodel is fully harnessed, for instance in terms of access to statistics of interest.

\subsubsection{Pseudospectral polynomial surrogate}
In this section, we consider without loss of generality that our quantity of interest is a scalar that depends on a single random parameter (e.g. $u(x,\xi)$).
Considering this quantity of interest  is sufficiently smooth, we rely on spectral expansions to approximate it. A finite set of basis functions is considered, and the solution is projected onto the subspace spanned by these functions. Thanks to a probabilistic framework, the hierarchical nature of the  representation provides a convenient functional format to perform subsequent tasks: sensitivity analysis, calibration, etc... Here, the standard projection-based gPC approximations $ \Pcal(u)$ is used:\\

\begin{equation}
u(\boldsymbol{x},\xi)\approx \Pcal_l(u) \equiv  \sum^{q_l}_{j=0} \BigE \left [ u(\boldsymbol{x},\xi) \cdot \phi_j(\xi) \right ] \phi_j(\xi)  =  \sum_{j=0}^{{q_l}} u_j(\boldsymbol{x}) \phi_j(\xi ),
\label{eq:projection}
\end{equation}
where $\phi_{j}$ are orthonormal Legendre polynomials and $ \Pcal_l$ is the polynomial truncation at level $l$. In higher dimensions $N$, its multivariate polynomial approximation is obtained via full tensor products. 
In its pseudospectral form, expansion coefficients will be calculated thanks to quadratures/cubatures that constitute a structured sampling of the parameters space. This comes at a certain cost but guarantees more accurate results than those obtained from empirically designed regression approaches \cite{Eck_2015}. In particular, the approximation will be subjected to truncation error but will ensure minimum aliasing error \cite{conrad}.

\subsubsection{Sparse sampling}
In this work, we attempt to minimize the number of deterministic simulations to perform. We rely on the use of sparse grids \cite{SG_book} in the framework of the pseudospectral polynomial approximation, which is different from the sparse grid stochastic collocation method of \cite{Chen2013}. This approach allows a careful selection of some representative sample points to form a sparse collocation in order to alleviate the computational effort and lower the cost of the prescribed high-dimensional integral computations \cite{Smolyak}.
Sparse grids that are structured in the form of quadratures/cubatures make error analysis more convenient because they integrate exactly functionals that are polynomials of a certain degree. They are incrementally assembled from a sequence of one-dimensional quadrature formulas. Many quadrature families exist with different (linear or exponential) growth of their nodes with respect to the quadrature level, and different polynomial integration capabilities. For this study, we rely on \textit{nested} sparse grids based on Gauss-Patterson-Kronrod (GP) \cite{Patterson1968} formulae that are known to be quite efficient.

\subsubsection{Moments and global sensitivity analysis}
Once the modal gPC coefficients are computed and with the help of the basis orthogonality, moments, confidence intervals, sensitivity analysis and probability density function of the solution can be readily evaluated \cite{olm_book2010}. \\ 
In particular, global variance-based sensitivity analysis may be performed in order to quantify, via correlation ratios,  the relative importance of each (or a group of) random input parameter to the uncertainty response of the system.
The Sobol' functional decomposition (or ANOVA decomposition) of $u=h(\bxi)$ is unique and hierarchic. We have:
\begin{equation}
u=h(\bxi)=\sum_{s\, \subseteq \{1,2,\ldots N_d\}} h_s (\bxi_s),
\label{eq:ANOVA}
\end{equation}
where $s$ is a set of integers such that $\bxi_s=(\bxi_{s_1},\ldots, \bxi_{s_N})$, with $N=\textrm{card}(s)=|s|$ and $h_{\varnothing}=h_0$. In this way, the variance of the solution can be derived from Eq. (\ref{eq:ANOVA}) and decomposed accordingly \cite{ANOVA_Efron81}:
\begin{equation}
\sigma^2=\sum_{s\, \subseteq \{1,2,\ldots N_d\}} \sigma_s^2, \quad \textrm{with} \quad \sigma_s^2 = \mathbb{V}\left ( \mathbb{E}\{ u | \bxi_s \} \right )- \sum_{ \substack{t \subset s \\ t \neq s \\ t \neq \varnothing} } \sigma_t^2,
\end{equation}
where $\mathbb{V}$ is the variance operator.
The normalized Sobol' indices \cite{Sobol_1993} $S_s$ are defined as:
\begin{equation}
\quad S_s \equiv  \frac{\sigma_s^2}{\sigma^2}  \quad \textrm{and}  \sum_{ \substack{  s\, \subseteq \{1,2,\ldots N_d\} \\ s \neq \varnothing}} S_s = 1,
\label{eq:sobol_coeff1}
\end{equation}
which measure the sensitivity of the variance of $u$ due to the interaction between the variables 
$\bxi_s$, without taking into account the effect of the variables $\bxi_t$  (for $t \subset s$ and $t \neq s$). First-order Sobol' indices were used in \cite{Chen2013,Eck_2015}. In this study, we compute first- and second-order indices.

\section{Numerical results}
\subsection{Validation of the hemodynamics deterministic solver}
Preliminary deterministic validations of our solver (not all reported here) have been performed against benchmarks  from the literature, in particular \cite{Formaggia_LCNSE2002,Boileau_2015}. 
The goal is to validate the solver for both piecewise-constant and spatially-varying arterial wall elastic properties.
\subsubsection{Aortic bifurcation with piecewise-constant elastic properties}
\begin{figure}
\centering
\includegraphics[width=9.5cm]{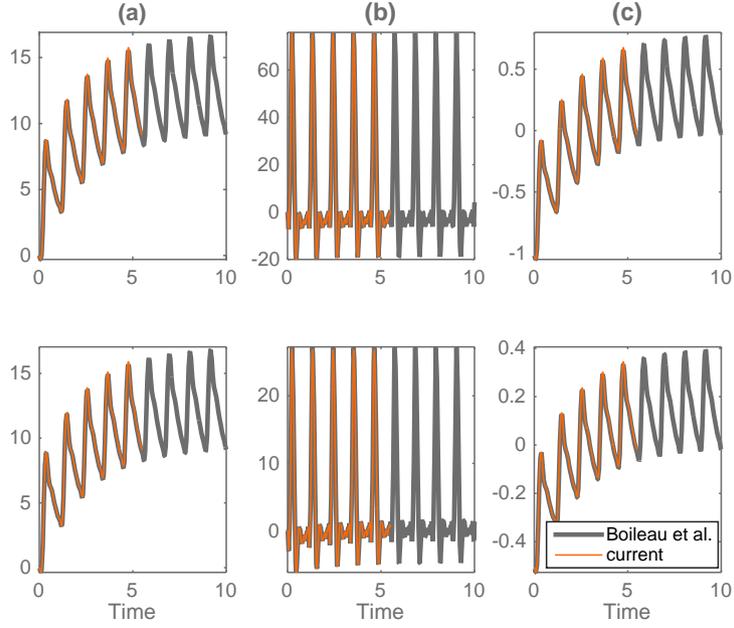}
\caption{Aortic bifurcation validation, cf. Boileau {\em et al.} \cite{Boileau_2015} Figure 5, p.17; (a): pressure $(kPa)$, (b): flow rate $(ml\cdot s^{-1})$ and (c): change in radius from diastole $(mm)$, at the midpoint
of the aorta (top) and midpoint of either iliac artery (bottom). 
}
 \label{fig:validation}
\end{figure}

The recent paper \cite{Boileau_2015} is a study of different numerical schemes for 1D blood flow modeling. Here, we focused more specifically on the comparison with the discontinuous Galerkin method for the case of intermediate complexity: namely the aortic bifurcation. This benchmark simulates the abdominal aorta branching into the two iliac arteries. It is a relevant computational experiment, that is easy to setup while retaining prominent network arterial features, i.e. -- a bifurcation model, -- {\em in vivo} inflow boundary condition signal, -- blood viscosity effect and -- 0D terminal boundary conditions with matched three-element Windkessel models. The results from our solver are compared in Figure~\ref{fig:validation}. It is clear that our solver performs adequately and does provide the expected solution.
\subsubsection{Vascular prosthesis in a single artery with variable elastic properties}
Reference cases are this time taken from this paper \cite{Formaggia_LCNSE2002}. The goal is to assess, in a single straight vessel of length $l_0$ with constant cross-section at rest and Young's modulus $E_0$, the effect of the local changes in vessel wall elastic characteristic due to the positioning of a prosthesis -- with $E_{0_{\mathrm{stent}}}=k\times E_0$ -- on the pressure pattern. These test cases are very severe because the stent is made of a much stiffer material than the surrounding artery. The variation in Young's modulus is not discontinuous but very steep. This causes a pressure build-up at the proximal point (P). Pressure histories at this point show that the interaction between the incoming and reflected waves reveals in discontinuities in the slope. All curves from the results of \cite{Formaggia_LCNSE2002} are successfully reproduced with our solver.
\begin{figure}
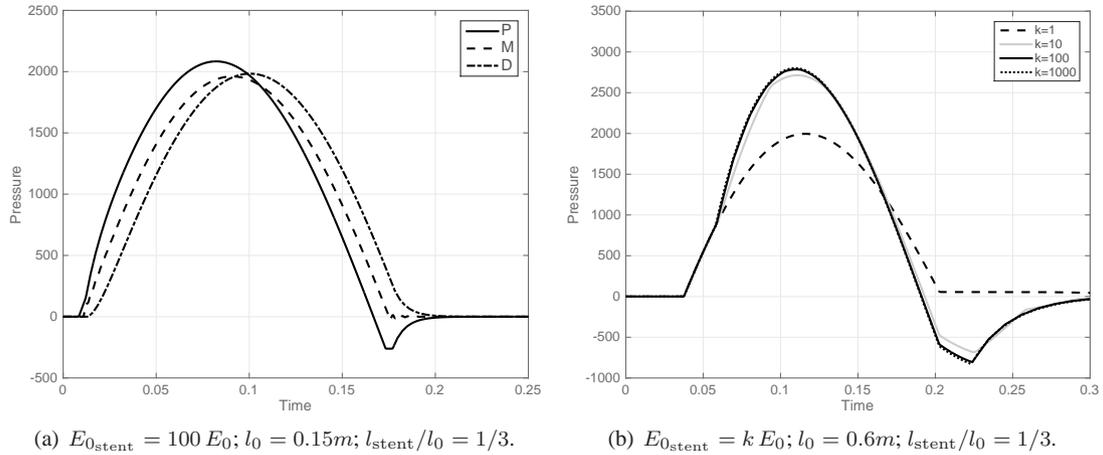

\centering
\subfigure[$E_{0_{\mathrm{stent}}}=100\, E_0$; $l_0=0.15 m$; $l_{\mathrm{stent}}/l_0=1/3$.]{\includegraphics[width=7.0cm]{./brault_figure7a.eps}}\hspace{0.1in}
\subfigure[$E_{0_{\mathrm{stent}}}=k\, E_0$; $l_0=0.6 m$; $l_{\mathrm{stent}}/l_0=1/3$.]{\includegraphics[width=7.0cm]{./brault_figure7b.eps}}
\caption{Vascular prosthesis validation: pressure history at different proximal (P), medium (M) and distal (D) monitoring points located prior, within and posterior to the stent position respectively, cf. Formaggia {\em et al.} \cite{Formaggia_LCNSE2002} Figures 7-8, p.151-152
.}
\label{fig:valid_stent}
\end{figure}

\subsection{Uncertainty quantification of the effect of stochastic inflow boundary condition on pulse waves response}\label{subsec:UQ_heart}

In the following, the influence of the inflow boundary conditions are presented  in terms of three main quantities, pulse pressure (PP),  arterial distensibility (AD) and  reflection magnitude (RM), which are respectively defined as:

$$
\left\{
\begin{aligned}
\textrm{PP}&=P_{\textrm{sys}}-P_{\textrm{dia}}\cr
\textrm{AD}&=(A_{\textrm{sys}}-A_{\textrm{dia}})/A_{\textrm{dia}}\cr
\textrm{RM}&= \textrm{P}_{\textrm{b}}/\textrm{P}_{\textrm{f}}\cr
\end{aligned}
\right.
$$

\noindent where  $\cdot_{\textrm{sys}}$ and $\cdot_{\textrm{dia}}$ refer to the systolic and diastolic values respectively and
 $\textrm{P}_{\textrm{f}}$ and $\textrm{P}_{\textrm{b}}$ are the maximal amplitudes of the separated forward and backward waves, cf. section \ref{subsubsec:wave_separation}. For this study, we rely on level $(l=5)$ \textit{nested} sparse grids based on Gauss-Patterson-Kronrod (GP) grid and Smolyak algorithm, totalizing 351 simulations. Multivariate ($N=3$) Legendre polynomial chaos of total order $P=3$, summing up to 20 terms, are sufficient to produce converged statistical results.

\subsubsection{Global statistics of PP, AD and RM}

\begin{figure}[ht]
\centering
\includegraphics[width=10cm]{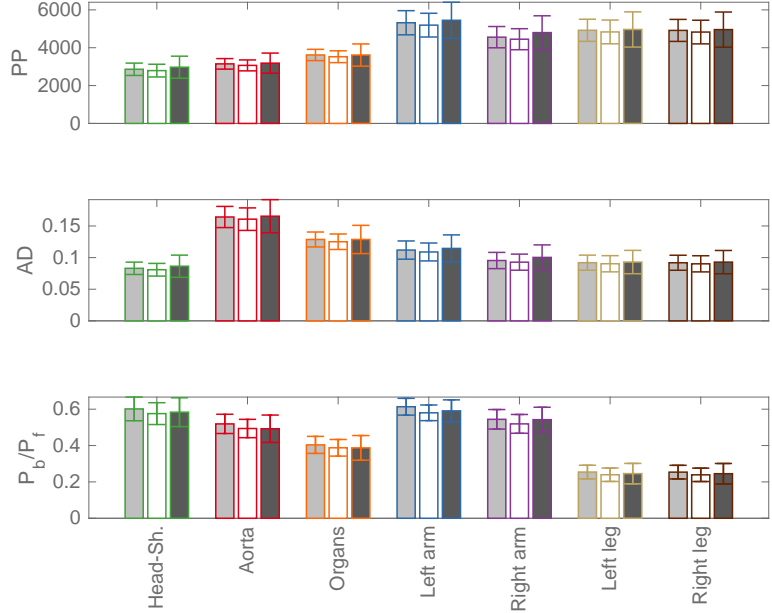}
\caption{Pulse pressure ($\mathrm{PP}_{[Pa]}$), area distensibility (AD) and reflection magnitude (RM$=\textrm{P}_{\textrm{b}}/\textrm{P}_{\textrm{f}}$) statistics (mean $\pm$ std) { averaged} by main cardiovascular groups. 
Light gray bars refer to the case of cardiovascular uncertainty modeling, while white (respectively dark gray) bars refer to the augmented case including proximal (respectively proximal-distal) aortic distensibility uncertainty modeling (see section \ref{subsec:UQ_heart+aorta}). \label{fig:PP_AD_stats}}
\end{figure}

Figure~\ref{fig:PP_AD_stats} shows the mean  statistics (mean and standard deviation) of PP, AD and RM averaged over various cardiovascular groups, with respect to inflow (and aortic stiffness) variations. Deterministic reference results at nominal condition are not included but departs from the mean solutions due to 
the nonlinear responses of the arterial tree. 

In particular, PP is a very complex quantity that depends on anthropometric factors such as age, gender, height and heart rate. It is determined not only by the ejection pattern of the left ventricle but also by aortic stiffness and the amplitude and timing of wave reflections. As expected, the results show that mean values of PP propagate along the arterial tree, with higher PP values in peripheral (e.g. brachial) arteries than in central (e.g. aortic and carotid) arteries.  
Standard deviations results show that PP uncertainty is always more pronounced  at distal locations. Averaged PP statistics are similar in the aorta and head/shoulders groups.

The next results from Figure~\ref{fig:PP_AD_stats} are relative to AD. As physiologically expected, they show that on average the most compliant part of the network is the aorta with quite a large variability due to its length. Moving to the periphery, the other groups are in average less and less compliant. We notice a slightly less significant variance of AD in the organs.\\
Finally, the last results from Figure~\ref{fig:PP_AD_stats} are relative to RM. This ratio is an accurate marker of the magnitude of wave reflections that are thought to be an important mechanism of augmentation of blood pressure with aging and in disease. The work of Baksi {\em et al.} \cite{JohnBaksi2016441} assesses the attenuation of reflected waves in man during retrograde propagation from the femoral artery to the proximal aorta. Twenty subjects underwent {\em invasive} measurement of pressure and flow velocity with a sensor-tipped intra-arterial wire at multiple locations distal to the proximal aorta before, during and following occlusion of the left femoral artery by thigh cuff inflation. 
Several of our findings are coherent with the results from this study: -- averaged values of RM close to the heart, e.g. at the proximal aorta, are close to their measurements, -- values of RM in the lower limbs are lower than aortic values which may be explained by re-reflections due to the marked impedance mismatch of bifurcations traversed in the retrograde direction. 

\subsubsection{Detailed statistics with respect to the cardiac output uncertain parameters}

\begin{figure}[ht]
\begin{center}
\includegraphics[width=7cm]{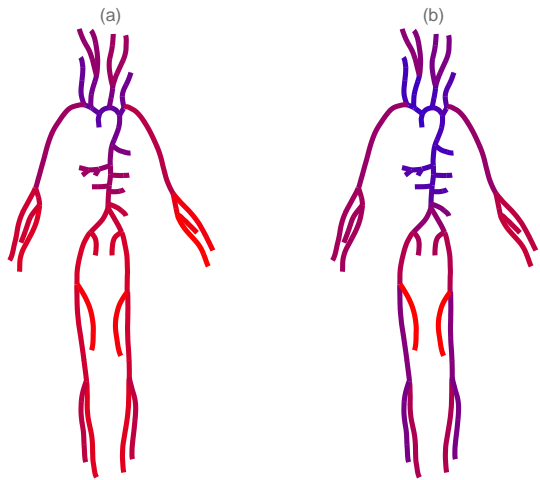}
\caption{Detailed pulse pressure PP vascular network  statistics: mean (a) and std (b). Low (respectively high) values are represented in blue (respectively red) color.  \label{fig:Network_PPmeanstd}
}
\end{center}
\end{figure}

\begin{figure}[ht]
\begin{center}
\includegraphics[width=10cm]{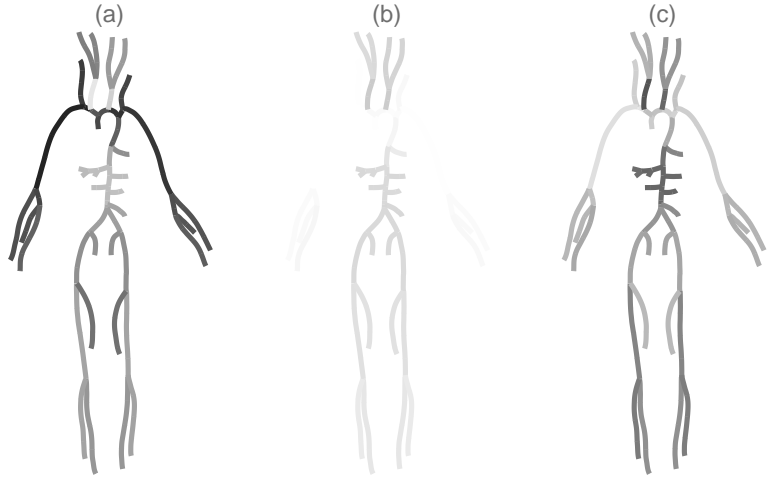}
\caption{Detailed PP first-order Sobol' indices: $S_{T}$ (a), $S_{\overline{Q}}$ (b) and $S_{Q_{\mathrm{max}}}$ (c). Low (respectively high) values are represented in white (respectively black) color.  \label{fig:Network_PPsobol}
}
\end{center}
\end{figure}

Due to the significant error bars at distal locations observed, it is necessary to refine our analysis. For instance, we wish to determine if PP amplification is  likely to be more pronounced than in nominal conditions and if it is more sensitive to some of the physiological parameters involved in the pulsatile cardiac output. Figures ~\ref{fig:Network_PPmeanstd} \&~\ref{fig:Network_PPsobol} display detailed distributions of PP low-order statistics and sensitivity analysis in the arterial network. Mean values of PP noticeably increases toward the periphery. We observe in particular a potent gradient along the aorta. Large mean values are obtained at the limb extremities. Std values are also larger in the limbs. However, statistical variations are quite low all along aorta and organs.\\
Sobol' indices show that the heart rate is very influential on the variability of the pulse pressure in the proximal region and the upper limbs, cf. Figure \ref{fig:Network_PPsobol}-(a). Maximum cardiac flow rate affects pulse pressure dispersion further down, in the lower part of the aorta, the organs and the lower limbs (c), while mean cardiac flow rate has little effect overall (b).\\

\begin{figure}[h]
\begin{center}
\includegraphics[width=10.0cm]{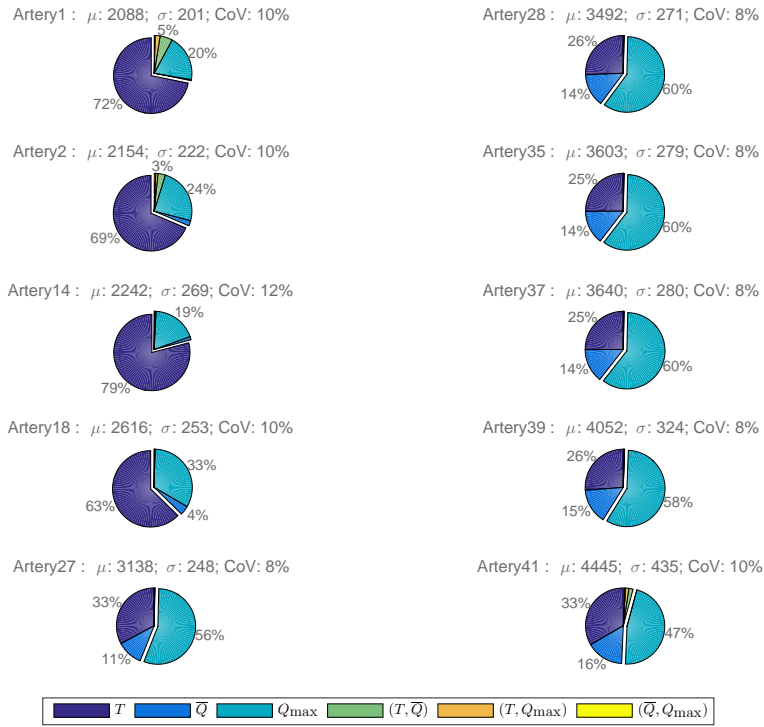}
\caption{Pulse pressure (PP) Sobol' indices distribution within the aorta arterial model segments \cite{Sherwin_2003a}.\label{fig:PP_aorta_pies}}
\end{center}
\end{figure}

\begin{figure}[ht]
\begin{center}
\includegraphics[width=5.5cm]{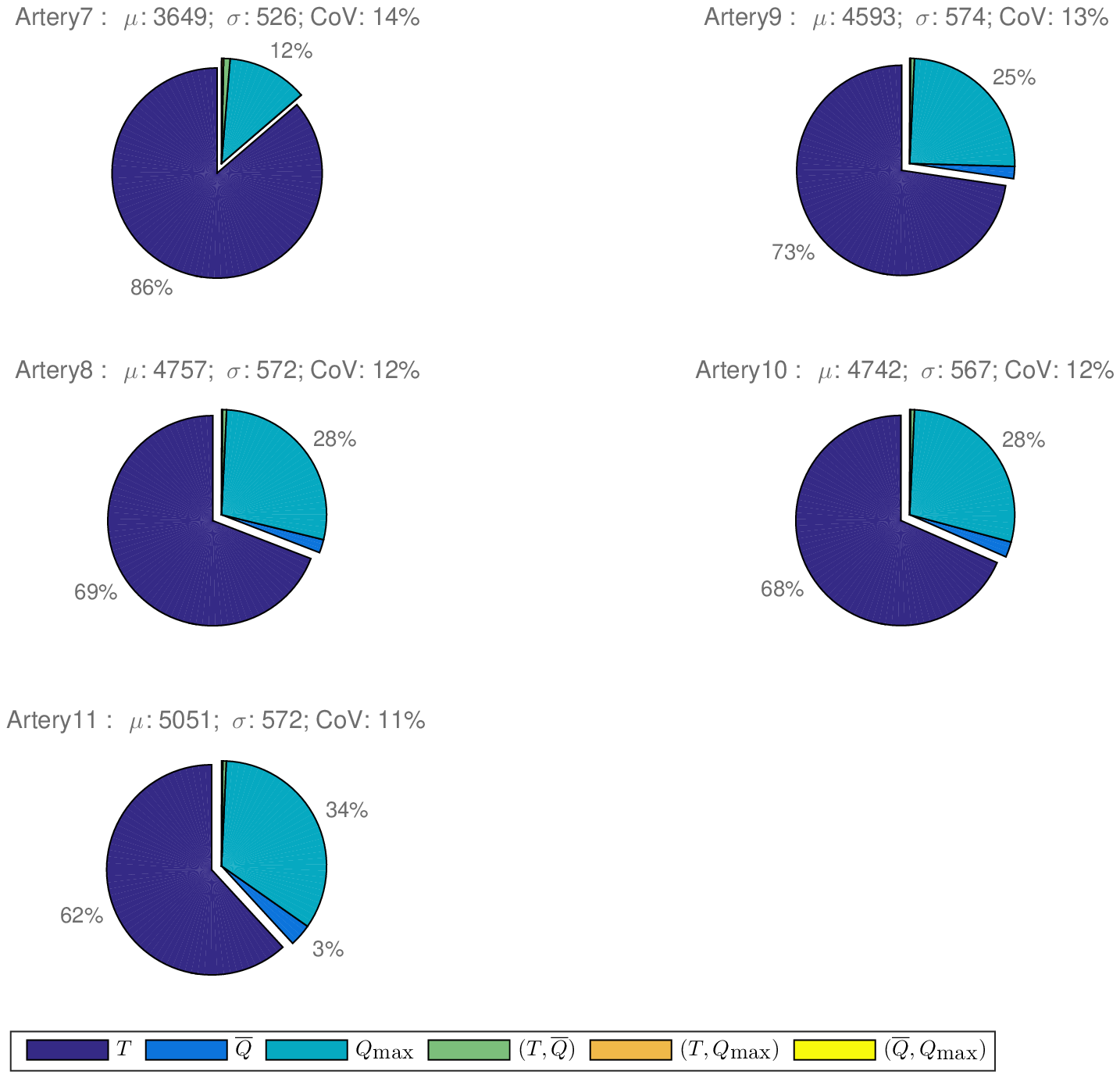}
\caption{Pulse pressure (PP) Sobol' indices distribution within the right upper limb model \cite{Sherwin_2003a}.\label{fig:PP_larm_pies}}
\end{center}
\end{figure}

\begin{figure}[ht]
\begin{center}
\includegraphics[width=8.0cm]{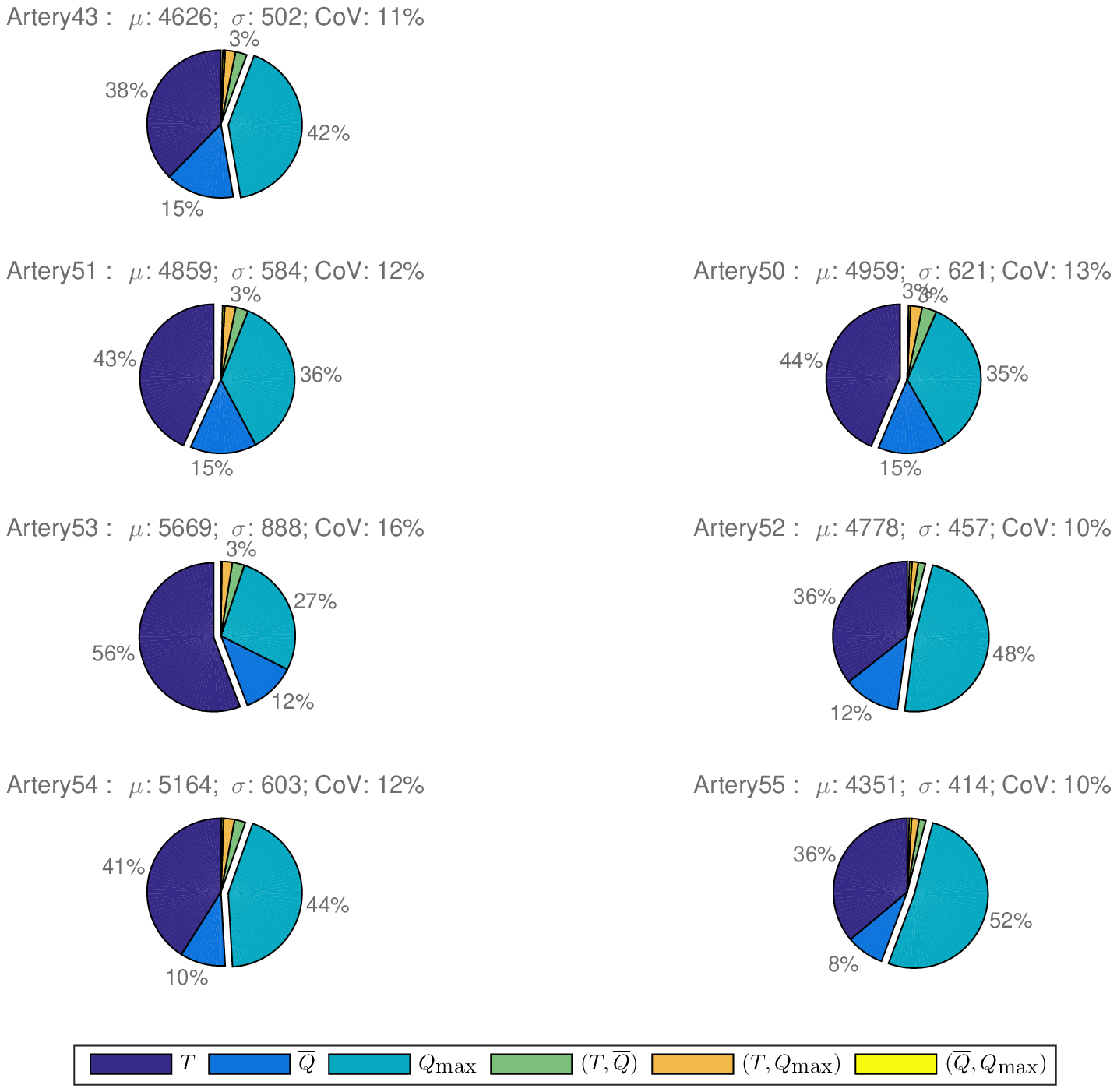}
\caption{Pulse pressure (PP) Sobol' indices distribution within the right lower limb model \cite{Sherwin_2003a}.\label{fig:PP_lleg_pies}}
\end{center}
\end{figure}

More quantitative results in some cardiovascular groups are more specifically examined in:  
 -- the aorta,  
 -- the (left) upper limb and  
 -- the (left) lower limb.
 The corresponding  figures~\ref{fig:PP_aorta_pies}--\ref{fig:PP_lleg_pies}  give the contribution of each kind of uncertainty, first or second order for the pulse pressure for each cardiovascular group.

\begin{figure}[ht]
\begin{center}
\includegraphics[width=9cm]{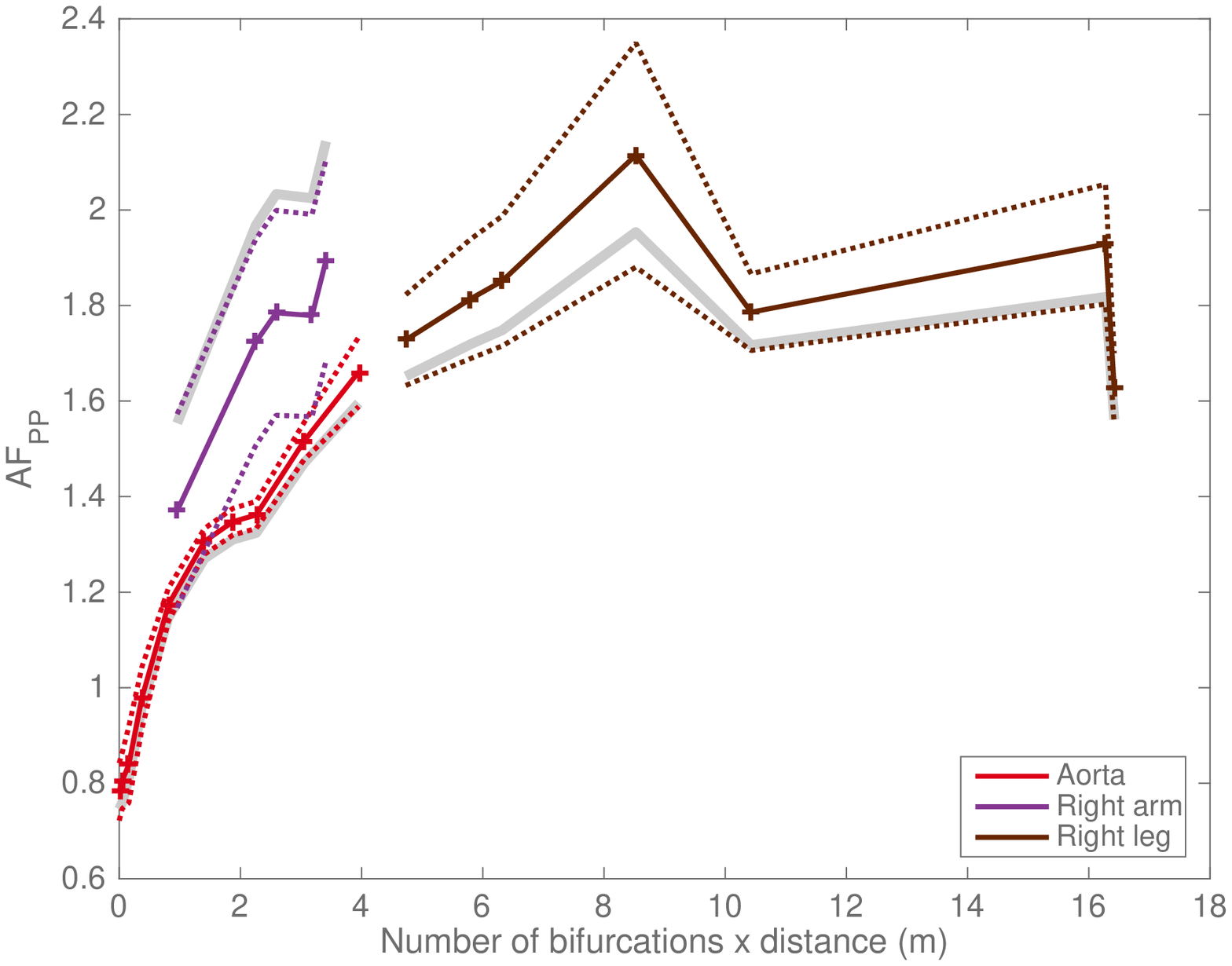}
\caption{Pulse pressure amplification factor AF$_{\mathrm{PP}}=\mathrm{PP}_{\mathrm{A}_{\{5\}}}/\mathrm{PP}_{\mathrm{A}_i}$ (i.e. relative to the right carotid) along the arterial network versus the number of bifurcations times the metric distance from the heart. Solid colored lines are the mean values and dashed colored lines are mean $\pm$ one std; solid gray lines are corresponding deterministic results at nominal conditions. 
\label{fig:amplification_factor}}
\end{center}
\end{figure}

In the upper part of aorta (i.e. upstream of the organs), the variability of the PP is mainly dominated by the effect of the heart beat uncertainty ($\sim 70\%$) and to some extent by $Q_{\textrm{max}}$ (less than $25\%$) whereas the effect of $\overline{Q}$ is negligible.  In the lower part of aorta (i.e. downstream of the organs), the variability of the PP is now mainly dominated by the effect of $Q_{\textrm{max}}$ ($\sim 60\%$) and to some extent by the heart beat uncertainty ($\sim 25\%$). The effect of $\overline{Q}$ is not negligible anymore but remains moderate close to $15\%$. Overall, for aortic PP, Sobol indices show that the different effects are mainly additive and there is little contribution from coupled interactions between input parameters. \\
The effect of cardiac variability on the distribution of peripheral and central pressure as well as pulse pressure amplification if of utmost importance for the selection of optimal measurement sites for better prediction of cardiovascular risk \cite{Wilkinson_JoP2000} . Figure~\ref{fig:amplification_factor} shows the evolution of the PP amplification factor (AF$_{\mathrm{PP}}$) statistics along some of the arterial network
paths. For each simulation, the quantity is computed relative to the pulse pressure in the (right) carotid midpoint, i.e. $\mathrm{PP}_{\mathrm{A}_{\{5\}}}$,  as it is often clinically done. The horizontal axis refers to the distance times the number of bifurcation levels from the aortic root.
The influence of the cardiac modeling uncertainty is potent when the results are compared to
 reference amplification factor at nominal conditions (gray solid lines). Deterministic AF$_{\mathrm{PP}}$ responses over predict the mean distributions in the upper limbs but under predict the mean distributions in the aorta and the lower limbs.
We notice the increase of the mean amplification  toward the network periphery, with an increase close to $+35\%$ in average for the brachial PP for instance. It reaches larger  values in the upper and lower limb extremities. Standard deviations are large all along the upper and lower limbs and grow moderately  away from the aorta root.

\subsubsection{Some statistical trends and correlation results}

\begin{figure}[ht]
\begin{center}
\includegraphics[width=10cm]{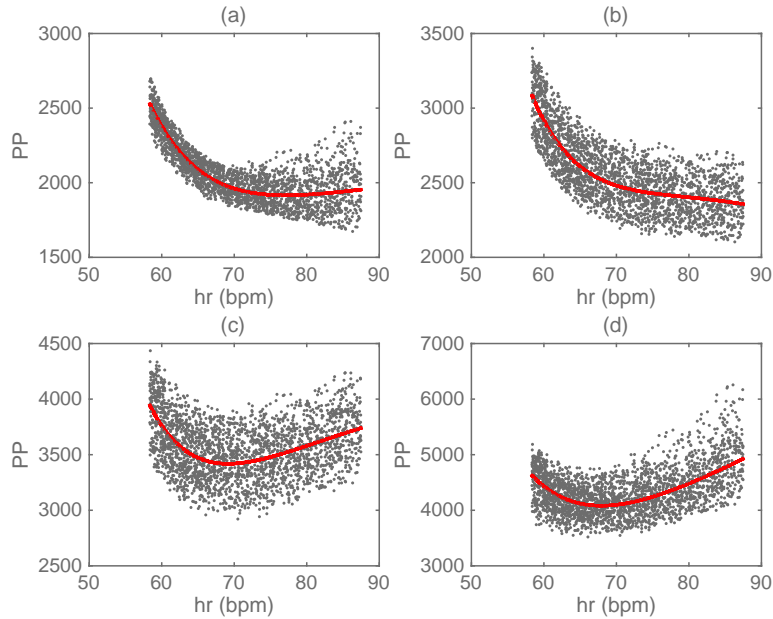}
\caption{Aortic PP dependance with respect to random variables $(hr,\overline{Q},Q_{\textrm{max}})$ (dotted samples: here projected onto the plane and solid red line: marginalized over $(\overline{Q},Q_{\textrm{max}})$) for (a): ascending aorta (A$_{[1]}$); (b): thoracic aorta I (A$_{[18]})$; abdominal III artery (A$_{[37]})$; abdominal V artery (A$_{[41]}$)). \label{fig:PPvsT_aorta}}
\end{center}
\end{figure}

The dependance of PP with respect to its main parameter in the aorta, namely $h_r$,  is presented on Figure \ref{fig:PPvsT_aorta}.

In the upper part of the aorta, pulse pressure has the tendency to decrease in a nonlinear fashion for increasing heart rates within the studied range, cf. plots (a)-(b) from Figure~\ref{fig:PPvsT_aorta}. Similar trends have been observed for { carotid} PP experimentally measured \cite{Albaladejo2001}.  In the lower part of the aorta, pulse pressure growth is not monotonic and reaches in average a minimum for $hr \sim 65-70$bpm, cf. plots (c)-(d) from Figure~\ref{fig:PPvsT_aorta}.

\begin{figure}[ht]
\begin{center}
\includegraphics[width=9cm]{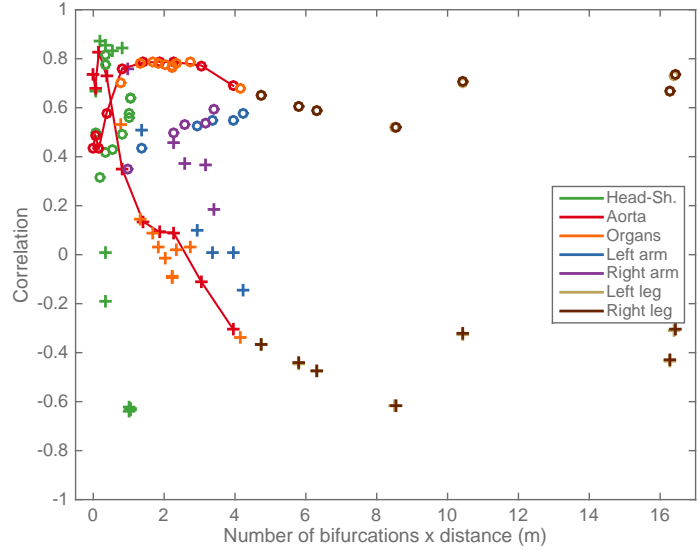}
\caption{Correlation coefficients $\rho_{\textrm{PP},T}$ ($+$) and $\rho_{\textrm{PP},Q_{\textrm{max}}}$ ($\circ$)  along the arterial network versus the number of bifurcations times the metric distance from the heart. Solid lines point out the aortic main arteries ordered from the proximal to distal alignment. \label{fig:correlationPP_T&Qmax}}
\end{center}
\end{figure}

The correlation coefficients between PP and its most influential inflow parameters such as $T$ and $Q_{\textrm{max}}$  are depicted on Figure \ref{fig:correlationPP_T&Qmax} and are also very revealing.
 In particular, it is clear that there is a good correlation of PP and $T$ in the upper part of the aorta -- cf. solid red lines in Figure~\ref{fig:correlationPP_T&Qmax} -- that drops considerably further down, while it is the opposite trend for the correlation with  $Q_{\textrm{max}}$. Interestingly PP in lower limbs and organs are much better correlated with $Q_{\textrm{max}}$ than $T$. In the head and shoulders region, brachiocephalic, subclavian and vertebral arteries are well correlated to $T$ while carotids are not as correlated and other arteries are even anticorrelated with $T$ such as internal/external carotids.

\subsection{Uncertainty quantification of inflow boundary condition and proximal arterial stiffness coupled effect on pulse waves response}
\label{subsec:UQ_heart+aorta}
Uncertainty modeling is now extended to the mechanical properties of the arterial walls in the proximal region of the heart. More specifically,  elastic characterization of the {\em aortic} vessel is modeled as a spatial random field with imposed low-order statistics and correlation, as explained in section \ref{subsubsec:aort_rand}. In this case, the uncertainty will affect the entire span of the aorta  -- represented in cyan color in figure~\ref{fig:network} -- i.e. from the proximal to the distal region adjacent to the main iliac bifurcation. We will refer to this case as the ``proximal-distal'' aortic distensibility uncertainty modeling case (or Case$_2$). For this study, the discretization being consequent on our modeling choice involves ($N=6$) random variables $\{T,\, \overline{Q},\, Q_{\mathrm{max}},\, c_{0_1},\, c_{0_2},\, c_{0_3} \}$, cf. Eq. (\ref{eq:KL_expansion}). For the sampling of this parametric space, 
we rely on level $(l=4)$ \textit{nested} sparse grids based on Gauss-Patterson-Kronrod (GP) rules and Smolyak algorithm, totalizing 545 simulations. A level $(l=5)$ sparse grid, i.e. 2561 simulations, was tested but did not modify significantly the results. Multivariate Legendre polynomial chaos of total order $P=3$, summing up to 84 terms, are required to produce converged statistical results.\\
Some UQ studies of arterial stiffening have suggested that uncertainty of the proximal aorta stiffness very close to the heart, maybe the leading cause of variability in the network pulse wave propagation response. Therefore, we will also consider a reference case where only a small region of the aorta near the system inlet boundary condition -- i.e. segment $\mathrm{A}_{\{1\}}$ in figure~\ref{fig:network}, corresponding to the distance materialized by the vertical black dashed thin line in Figure~\ref{fig:KL_decomposition}-(b) -- bears a random stiffness. We will refer to it as the ``proximal'' aortic distensibility uncertainty modeling case (or Case$_1$). For this study, with ($N=4$) random variables $\{T,\, \overline{Q},\, Q_{\mathrm{max}},c_{0_{\mathrm{A}_{\{1\}}}} \}$, we rely on level $(l=5)$ \textit{nested} sparse grids based on Gauss-Patterson-Kronrod (GP) grid and Smolyak algorithm, totalizing 769 simulations. Multivariate Legendre polynomial chaos of total order $P=3$, summing up to 35 terms, are sufficient to produce converged results. The $c_{0_{\mathrm{A}_{\{1\}}}}$ random variable is uniformly distributed with constant mean and std values that are chosen to match the statistics of the proximal-distal aortic distensibility distribution, averaged over the $\mathrm{A}_{\{1\}}$ span.

\subsubsection{Global statistics of PP, AD and RM}
A quick look at the statistics of Figure~\ref{fig:PP_AD_stats} is enough to demonstrate that the main trends are conserved when extending uncertainty quantification to the effects of aortic stiffness. For each cardiovascular group, Case$_1$ and Case$_2$ depart from the reference case described in section \ref{subsec:UQ_heart}. In particular, mean and std results are always noticeably larger for PP and AD in Case$_2$. This statement indicates that the increased level of uncertainty in the system, due to the aortic stiffness modeling results in larger amplitude variabilities of pressure and arterial distensibility. This is particularly true in peripheral regions. Larger std distributions of reflection magnitude RM are also observed in Case$_2$ but the mean values are lower than the reference case. This is interesting as it seems that the spatial variation of pulse wave velocities, due to the variations in arterial stiffnesses along the aorta, affects the magnitude of the reflected pressure waves relative to the transmitted ones.  

\subsubsection{Sensitivity analysis of the waves reflection magnitude}
\begin{figure}[ht]
\begin{center}
\includegraphics[width=10.0cm]{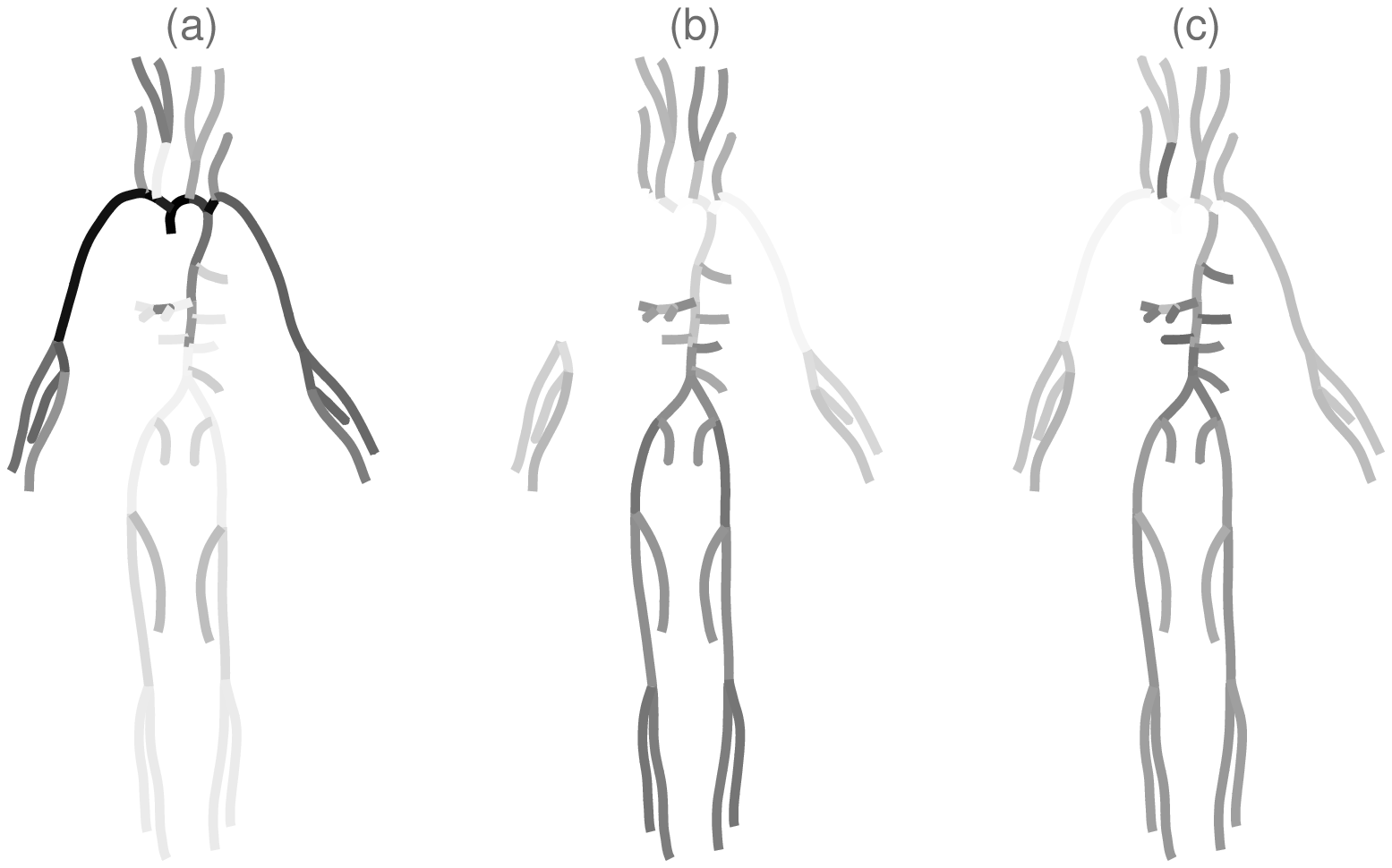}\caption{Detailed RM first-order Sobol' indices: $S_{T}$ (a), $S_{\overline{Q}}$ (b) and $S_{Q_{\mathrm{max}}}$ (c), reference case.\label{fig:Network_RMsobol1}}
\includegraphics[width=10.0cm]{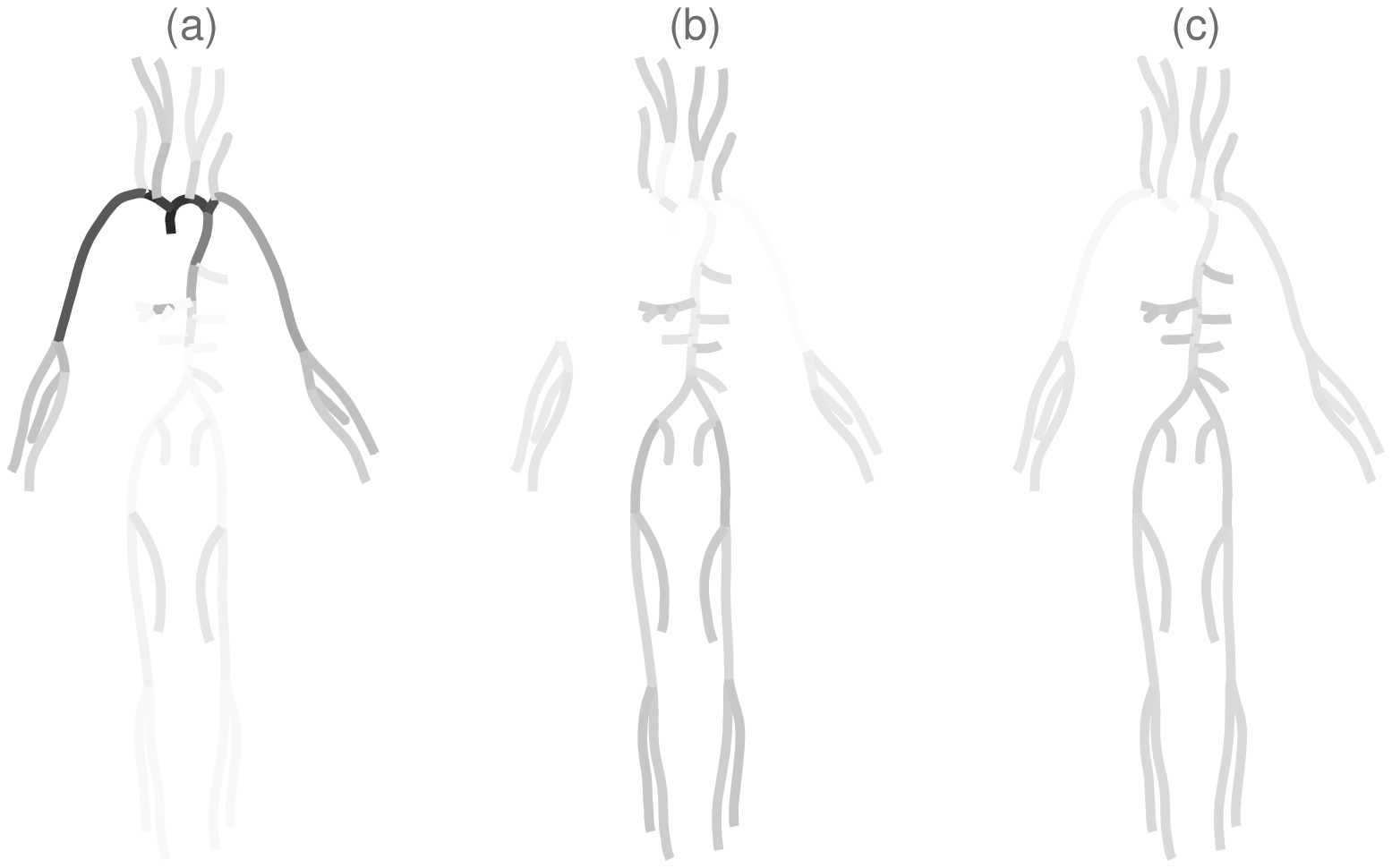}\\
\includegraphics[width=10.0cm]{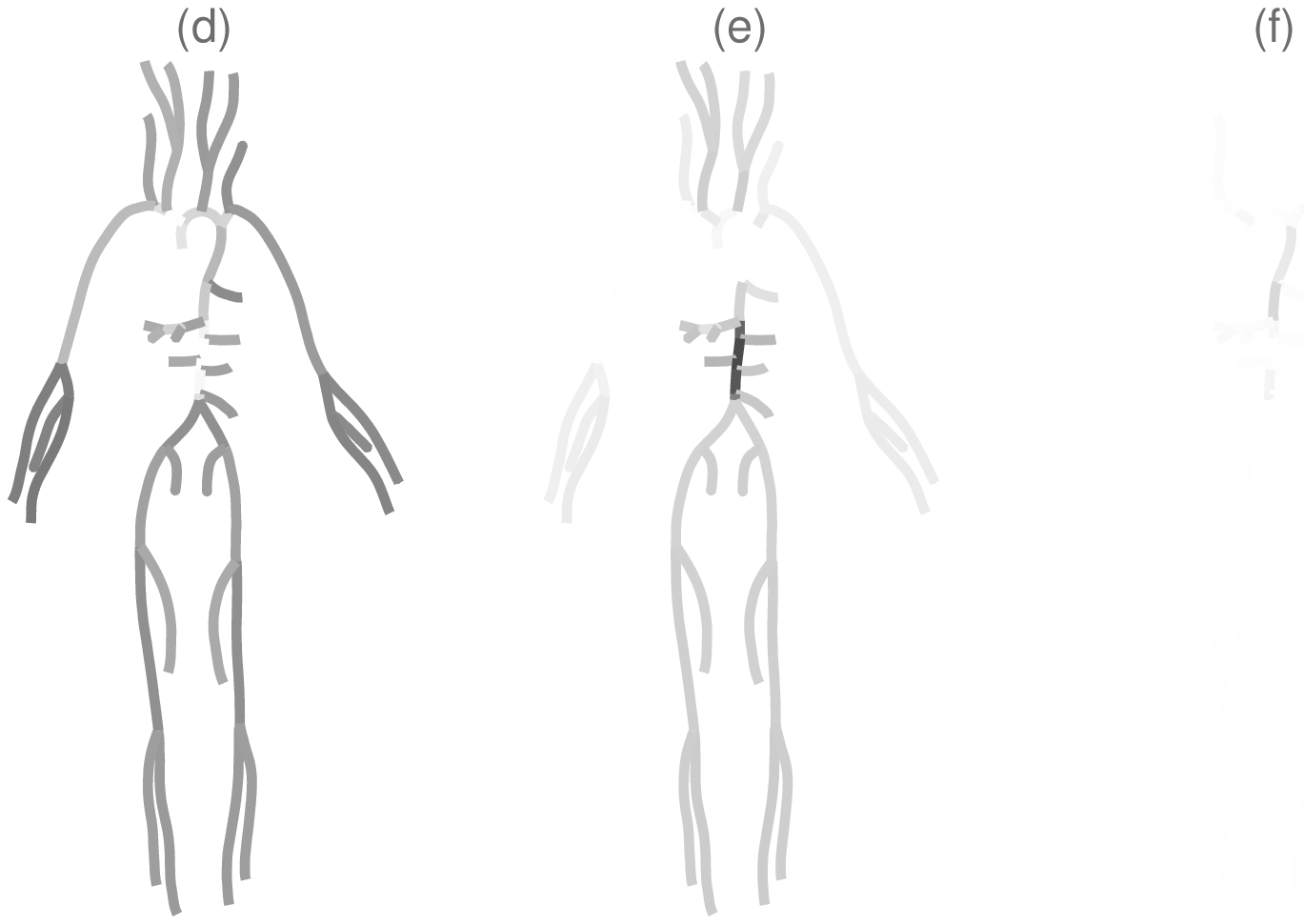}\caption{Detailed RM first-order Sobol' indices: $S_{T}$ (a), $S_{\overline{Q}}$ (b), $S_{Q_{\mathrm{max}}}$ (c), $S_{c_{0_1}}$ (d), $S_{c_{0_2}}$ (e) and  $S_{c_{0_3}}$ (f) for Case$_2$. Low (respectively high) values are represented in white (respectively black) color.  \label{fig:Network_RMsobol2}}
\end{center}
\end{figure}

A natural concern is the one of the spatial receptivity of the augmented system with respect to the cardiac variability: do we notice some modification in the reflection sites sensitivity to cardiac parameters when aortic stiffness uncertainty is taken into account? Figure~\ref{fig:Network_RMsobol1} and Figure~\ref{fig:Network_RMsobol2} compare the RM first-order Sobol' coefficients of the reference case of section \ref{subsec:UQ_heart} with the results of Case$_2$ respectively. In the reference case, largest values of RM std are obtained close to the heart at the aortic root, and in the arms to some extent (results not included). The std magnitude globally decreases toward periphery. Locations of large std are dominated by the parameter $T$. In other regions, influence of both $\overline{Q}$ and $Q_{max}$, with a slight advantage to the last one in aorta lower part and organs. \\
The injunction of spatial uncertainty in the modeling of the aortic wall stiffness is very significant as it strongly impedes on the sensitivity of the reflections within and outside the aorta region. In Figure~\ref{fig:Network_RMsobol2}, we notice that the sensitivity to $\overline{Q}$ and $Q_{\mathrm{max}}$ is flushed out from the system and reattributed to the random variables that describe the modal decomposition of the aorta distensibility, i.e. the $c_{0_i}$ modes. We recall that these random variables emerge from the discretization of the random field and are hierarchically distributed to form its low-dimensional model. They  correspond to the random magnitudes of the characteristic spatial stiffness fluctuations/scales  identified though their eigenspectrum. 
We see that the reflections become sensitive to the variation of the first mode in the peripheral regions: i.e. forearms, lower limbs and to some extent in the carotids, while the sensitivity of the aorta to this scale is very low (d). Interestingly, the second mode has a very predominant and local effect along the lower part of the aorta (e) despite the linearly decreasing values of stiffness variability imposed in this region, cf. Figure~\ref{fig:KL_decomposition}-(b). The amplitude of the third mode is too low to affect the reflections which are insensitive to this component (f).\\ Even if the physical interpretation of these distributions is difficult, it is interesting to observe how the pressure wave reflections vary in response to the aortic stiffness random fluctuations: the sensitivity to cardiac features is reallocated across the network and only the heart rate plays a dominant role in the proximal aorta and right upper limb and shoulder. The first KL eigenmode, with the largest wavelength, affects the entire tree but the aorta, while the second eigenmode with shorter wavelength alters the reflections in the aorta lower portion exclusively. 

\subsubsection{Correlation results}
\begin{figure}[ht]
\begin{center}
\includegraphics[width=9cm]{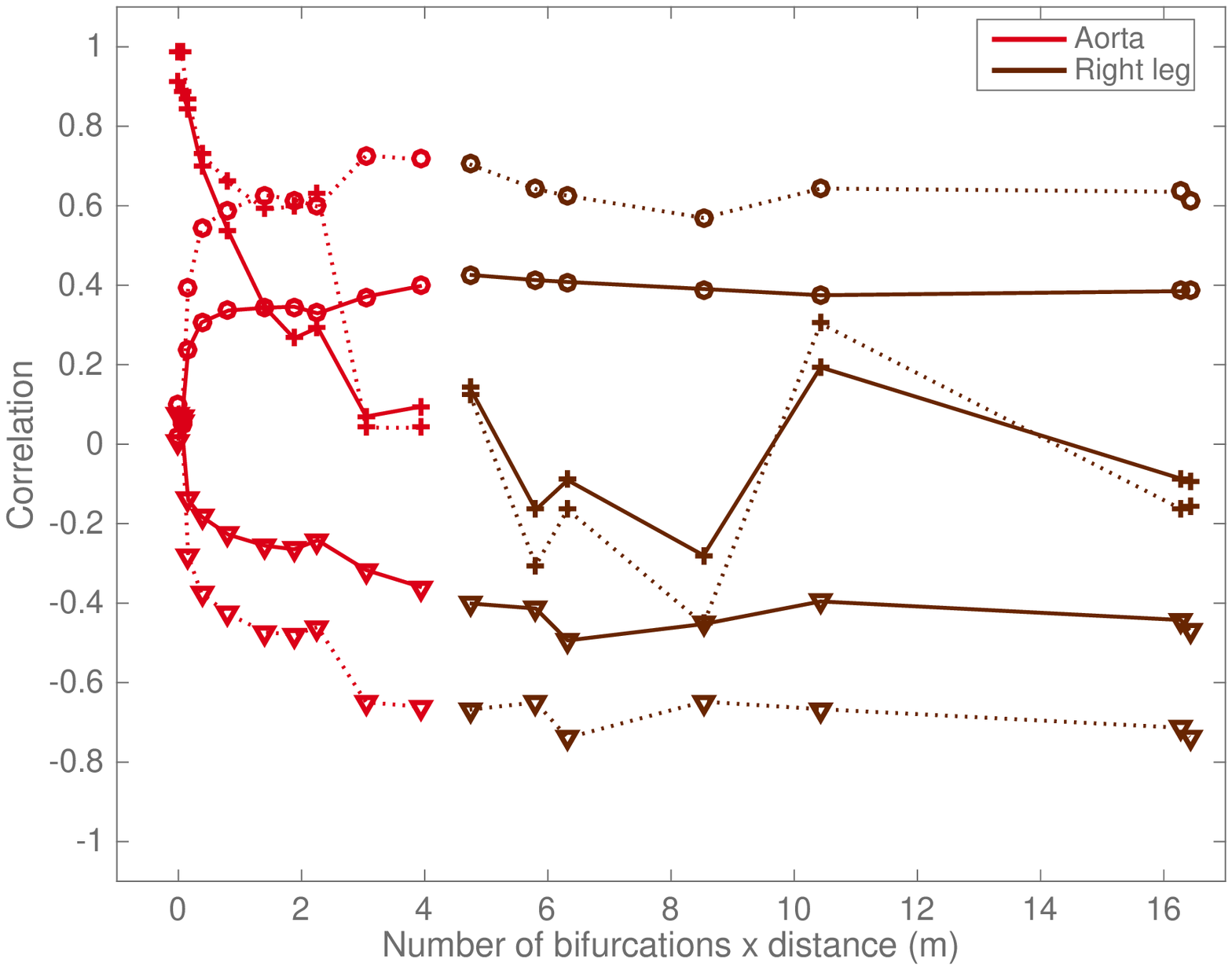}
\caption{Correlation coefficients $\rho_{\mathrm{RM},T}$ ($+$), $\rho_{\mathrm{RM},\overline{Q}}$  ($\triangledown$) and $\rho_{\mathrm{RM},\,Q_{\textrm{max}}}$ ($\circ$)  along part of the arterial network. Dotted lines refer to the case of cardiovascular uncertainty modeling alone, while solid lines refer to the augmented case accounting for proximal-distal aortic distensibility uncertainty modeling. \label{fig:correlationPbPf_T&Qmax}}
\end{center}
\end{figure}
The distributions of correlation coefficients between the reflection magnitude RM and the cardiac uncertain parameters are depicted on figure \ref{fig:correlationPbPf_T&Qmax}, along the aorta and the right lower limb. A strong positive correlation between RM and the cardiac cycle period $T$ exists in the proximal aorta. For the case of deterministic aortic distensibility, this correlation mildly drops along the aorta span while positive (respectively negative) correlation builds up with $Q_{\textrm{max}}$ (respectively $\overline{Q}$). These levels of correlation at distal location in the aorta remain sensibly  identical further down in the legs. Interestingly, in the case of aortic random spatial distensibility, correlation coefficients are qualitatively distributed in a similar manner but magnitudes are strongly attenuated. Correlations between RM and heart rate, decay for instance much faster along the aorta. Correlations between RM and $Q_{\textrm{max}}$ and $\overline{Q}$
are also less significant.

\section{Discussion and conclusion} 
In this study, we have relied on a stochastic sparse pseudospectral representation to quantify and analyse the effect of inherent uncertainties from cardiac output on the sensitivity of a human compliant arterial network response based on a reduced-order pulse wave propagation model. More specifically we have mainly focused on the sensitivity of two indices: -- the pulse pressure and -- the pulse reflection magnitude, which are known to be determined in part by the ejection pattern of the left ventricle, but also by aortic stiffness and other parameters. We were interested by the relationship between central and peripheral pulse pressure 
and the effect of wave reflection in the arterial system that is also well recognized although some aspects of the topic remain controversial.\\
In this work, we have not used a heart model but a simple pulsatile inlet flow waveform reproducing the most relevant cardiac features with a minimum number of  parameters: namely  cycle period, mean cardiac output and peak-to-mean flow ratio. We have also modeled another critical source of uncertainty arising from  the spatial heterogeneity of the aortic compliance. Indeed, stiffness properties of the aortic walls are not well known but play a key role in the propagation and damping of pulse waves generated at each cardiac cycle. A continuous representation of the aortic stiffness in the form of generic random fields of prescribed spatial correlation has been considered.\\
We have first validated our deterministic solver against specific benchmarks from the literature for both piecewise-constant and variable arterial walls stiffness. This step allowed us to be confident in our control of the deterministic error levels due to discretizations. 
We have then carried out simulations of the system under stochastic conditions.
Statistical results distribution across the network show that mean values are in agreement with physiological principles and clinical evidences while std magnitudes are of the same order of magnitude as the ones of the random input parameters.
The effects of the parametric uncertainty are mainly additive for the observables under consideration, but they reflect the nonlinearity of the model.
Mean cardiac output variability has a small and moderate impact on the pulse pressure and pressure waves reflection magnitude respectively. Heart rate variability, in the considered range, has a strong influence in the proximal region of the aorta, upper limbs and part of the head-shoulders cardio vascular group. Variability of the  flow rate maximum amplitude at the aortic root bears no more than a mild effect in the distal locations of the network and in the lower aorta portion.\\
Overall, the response in the aorta is very rich with: -- marked spatial gradients  of  both pulse pressure and reflection magnitude mean values, -- a negative  gradient of the reflection magnitude std values. 
Under deterministic physiological conditions, we expect that forward travelling waves from the aorta are well-matched at each bifurcation. This usually implies that reflected waves are not well-matched when propagating in the reverse direction toward the root of the arterial tree and this leads to the wave-trapping phenomena. Under stochastic conditions, the results suggest that the reflection magnitude in the lower generations of the bifurcating network is less sensitive to parameters inducing a change in wave propagation timing (such as heart rate) than the upper generations (close to the heart) due to the wave-trapping phenomena.\\
For the pulse pressure, std values are very stable along the aorta but maybe attributed to different parameters variability: i.e. cardiac beating period in the proximal region and maximum flow rate in the distal region respectively. This shift in the response sensitivity is well supported when looking at the statistical correlations. This balance is quite similar for the reflection magnitude.\\
The introduction of uncertainty in the structural properties of the aorta strongly modifies the statistical response, 
especially for a random {\em field} modeling that distributes stochastic fluctuations in a continuous manner all along its span. The consequence is a global increase of system response level of variability with a network remapping of the parametric sensitivity. The influence of the cardiac parameters noticeably decreases, due to declining correlations, except for the heart rate in the proximal aorta and right upper limb. For the reflection magnitude, the effects of aortic rigidity fluctuations (with long wavelength range) become important all across the system and not only downstream of the aorta, with an emphasis close to the body extremities. Again, the response along the aorta is very diverse and variable with a clear change from heart rate hegemony in the proximal region to aortic rigidity control (with medium wavelength range) in the distal region. The aortic uncertainty also induces changing stiffness gradients that directly affect the magnitude of wave reflections. It seems that the effects of the fluctuations with long wavelength range are only felt ``far away'' from their spatial origin. Fluctuations of shorter wavelength have a more local and direct effect on the reflections in a portion of the aorta itself. However, maybe due to their lower amplitudes or some complex wave-trapping phenomena, this effect vanishes in more distal locations. These last results will certainly deserve more fundamental and in depth investigations. \\
On a general viewpoint, deterministic and stochastic results presented here, hint at a noticeable dissymmetry between quantities predicted in the left and right arms. This is not the case for the lower limbs for which the statistics are very similar. To our knowledge, this difference was never reported in the literature relying on the use of this 55-segment arterial network topology. A careful analysis of our results suggests that this dissimilarity originates from the generation shift of the bifurcating network. Indeed blood ejected from the heart flows through four bifurcations to reach the upper left arm, but only three to reach the upper right one. \\
There are several limitations in our work: 1-- we have not used a model to describe the cardiac ejection dynamics 2-- the random parameters are assumed to be independent and to follow standard non-informative probability distribution functions,  3-- uncertainty modeling is incorporated in the description of the arterial wall stiffness of the aorta alone, 4-- this model necessitates a correlation form and typical length which are arbitrarily chosen. Items 2 and 4 could be improved from a careful statistical analysis of the clinical data and measurement campaigns reported in the literature, and could also benefit from the incorporation of a more complex and realistic heart model for Item 1. Item 3 also points to the lack of uncertainty modeling of the aorta cross-section, that would be required if one wants to account for the continuous wave reflections phenomena induced by the natural tapering of the vessel. In this case, additional constraints may arise due to the random fluctuations of the structural properties/geometry leading to non well-matched arterial bifurcations.
It remains that under stochastic conditions, physical interpretation and global understanding are intricate due to the complexity of the chosen closed distribution network and the fast stochastic interaction of different long-wave propagation and reflection features, that render individual wave tracking extremely strenuous.\\
Despite these limitations, we believe the results presented in this study have potential physiological and pathological implications. They will hopefully provide some guidance in clinical data acquisition and future coupling of arterial pulse wave propagation reduced-order model with more complex beating heart models. Finally, this work may also be seen as a contribution to the question of the relevance of generalized (or modeled) versus patient-specific (or measured) inflow boundary conditions of reduced-order pulse wave propagation model and more generally in multi-scale arterial hemodynamics coupling.

\smallskip
\acks{
The authors are greatly thankful to Dr J. Alastruey, from King's College London, for the fruitful exchanges related to some of the numerical validation aspects of this work.}

\section*{Conflict of interest statement}
Nothing to disclose.

\begin{table}[!ht]
\footnotesize
\begin{center}
\begin{tabular}{ | l | c | c | c |}
\hline  & $T$  & $\overline{Q}$  & $Q_{\textrm{max}}$    \\
\hline
\hline mean $\mu$  & 0.86 (s) & 100 (ml/s) & 650 (ml/s)     \\
\hline std $\sigma_{ [ \% \mu]}$ & $11.55\%$  & $5.77\%$  & $5.77\%$     \\
\hline
\end{tabular}
\end{center}
\caption{Inflow random variables statistics.}
\label{tab:stats}
\end{table}
\end{document}